\documentclass{article}

\usepackage{graphicx}

\usepackage{times}
\usepackage{mathptmx}
\usepackage{latexsym}
\usepackage{amsmath} 
\usepackage{amssymb}  

\usepackage{dsfont} 

\usepackage{url}  
\usepackage{epsf}
\usepackage{xspace}	
\usepackage{pgf}
\usepackage[lined]{algorithm2e}

\usepackage{amsthm} 
\newtheorem{theorem}{Theorem}
\newtheorem{proposition}[theorem]{Proposition}
\newtheorem{corollary}[theorem]{Corollary}
\newtheorem{lemma}[theorem]{Lemma}
\newtheorem{example}{Example}
\newtheorem{definition}{Definition}

\usepackage{amsopn}

\definecolor{red}{rgb}{0.8,0.2,0.2}
\definecolor{blue}{rgb}{0,0,0.5}
\definecolor{green}{rgb}{0,0.7,0}
\definecolor{violet}{rgb}{0.5,0.2,0.5}
\definecolor{orange}{rgb}{0.8,0.5,0.2}
\definecolor{meu}{rgb}{0.6,0.8,0.6}

\usepackage[colorlinks,pagebackref,linktocpage]{hyperref}

\newsavebox{\fmbox}

\newcommand{\xaxis}{$x$-axis\xspace}
\newcommand{\yaxis}{$y$-axis\xspace}
\newcommand{\auc}{\ensuremath{\mathit{AUC}}\xspace}

\newcommand{\mae}{\ensuremath{\mathit{MAE}}\xspace}
\newcommand{\mse}{\ensuremath{\mathit{MSE}}\xspace}

\newcommand{\OVER}{\ensuremath{\mathit{OVER}}\xspace}
\newcommand{\UNDER}{\ensuremath{\mathit{UNDER}}\xspace}
\newcommand{\RROC}{RROC\xspace}
\newcommand{\aoc}{\ensuremath{\mathit{AOC}}\xspace}
\newcommand{\mmse}{\ensuremath{\mathit{MMSE}}\xspace}

\newcommand{\absloss}{\ensuremath{{\ell}^A}}          
\newcommand{\squloss}{\ensuremath{{\ell}^S}}          

\newcommand{\aabsloss}{\ensuremath{{\ell}^A_{\alpha}}}   
   
\long\def\comment#1{}

\newcommand{\arxiv}[1]{{{#1}}}
\newcommand{\pr}[1]{{{}}}

\newcommand{\vect}[1]{{\bf{#1}}}

\begin{document}

\arxiv{\title{A graphical analysis of cost-sensitive regression problems}}
\pr{\title{ROC Curves for Regression}}  
\author{Jos\'e Hern\'andez-Orallo\hfill(jorallo@dsic.upv.es)\\
Departament de Sistemes Inform\`atics i Computaci\'o\\
Universitat Polit\`ecnica de Val\`encia, Spain\\}			  


\date{\today}


\maketitle

\begin{abstract}

\pr{Receiver Operating Characteristic (ROC) analysis is one of the most popular tools for the visual assessment and understanding of classifier performance. By means of ROC analysis we can analyse many issues such as performance for different operating conditions, regions where each classifier dominates, derive optimal thresholds and calculate aggregated measures, such as the Area Under the Curve (AUC), which is a good indicator of overall classifier performance.}
Several efforts have been done to bring ROC analysis beyond (binary) classification, especially in regression. However, the mapping and possibilities of these proposals do not correspond to what we expect from the analysis of operating conditions, dominance, hybrid methods, etc. In this paper we present a new representation of regression models in the so-called regression ROC (\RROC) space. The basic idea is to represent over-estimation on the \xaxis and under-estimation on the \yaxis. The curves are just drawn by adjusting a {\em shift}, a constant that is added (or subtracted) to the predictions, and plays a similar role as a threshold in classification. From here, we develop the notions of optimal operating condition, convexity, dominance, and explore several evaluation metrics that can be shown graphically, such as the area over the \RROC curve (\aoc). In particular, we show a novel and significant result, the \aoc is equal to the error variance (multiplied by a factor which does not depend on the model). The derivation of \RROC curves with non-constant shifts and soft regression models, and the relation with cost plots is also discussed.


\end{abstract}

\newpage
\tableofcontents
\newpage

\section{Motivation}\label{sec:intro}

\pr{ROC analysis is ... \cite{lusted1971,Swe86,Bradley1997, SDM00,flach2003decision,lasko2005use,Fawcett06,krzanowski2009roc} (radiology, medicine, statistics, bioinformatics, machine learning, pattern recognition, to name a few)}
In classification, the traditional notion of operating condition is common and well understood. Classifiers may be trained for one cost proportion and class distribution (both making the operating condition) and then deployed on a different operating condition. Some of the techniques and notions for addressing these cases are cost matrices, cost-sensitive classification \cite{Elk01} and very especially ROC analysis \cite{lusted1971,Swe86,Bradley1997, SDM00,flach2003decision,lasko2005use,Fawcett06,krzanowski2009roc}. 
ROC space decomposes the performance of a classifier in a dual way. On the \xaxis we show the false positive rate (FPR) and on the the \yaxis we show the true positive rate (TPR). ROC curves neatly visualise how the TPR and the FPR change for different (crisp) classifiers or evolve for the same (soft) classifier (or ranker) for a range of thresholds. The notion of threshold is the fundamental idea to adapt a soft classifier to an operating condition. ROC analysis is the tool that illustrates (among other things) how classifiers and threshold choices perform. The number and variety of applications and areas (radiology, medicine, statistics, bioinformatics, machine learning, pattern recognition, to name a few) have been increasing over the years \cite{Goin1982,Mamitsuka2006,Khreich2010,Khreich2012}. Also, some metrics derived from the ROC curve, such as the Area Under the ROC Curve (AUC) are now key for the evaluation and construction of classifiers \cite{ferri2002learning,Marrocco2008,Toh2008,hand2010evaluating,Ricamato2011,Kim2012}

The adaptation of ROC analysis for regression has been attempted on many occasions. However, there is no such a thing as the `canonical' adaptation of ROC analysis in regression, since regression and classification are different tasks, and the notion of operating condition may be completely different. In fact, the mere extension of ROC analysis to more than two classes has always been difficult because the degrees of freedom grow quadratically with the number of classes (see, e.g., \cite{srinivasan1999note,ferri2003volume,Schubert2011}).  
\arxiv{The inclusion of probabilities (and other magnitudes) in ROC curves or the use for abstaining classifiers \cite{ferri2004cautious,Pietraszek2005,ferri2005modifying} has not paved the way on how to do similar things for regression.} 
Consequently it is even questionable whether a similar graphical representation of ROC curves in regression (or other tasks \cite{rocai2004}) can even be figured out.
Notable efforts towards ROC curves (or graphical tools) for regression are the 
 Regression Error Curves (REC) \cite{bij2003regression}, the Regression Error Characteristic Surfaces (RECS) \cite{torgo2005regression}, the notion of utility-based regression \cite{torgo2009precision} and the definition of ranking measures \cite{rosset2007ranking}. These approaches are based on gauging the tolerance, rejection rules or confidence levels. Some of these approaches actually convert a regression problem into a classification problem (tolerable estimation vs. intolerable estimation).
Another recent approach has been based on the calculation of Kendall's rank $\tau$ correlation coefficient between the predicted and actual values \cite{ROCregrankers}, so disregarding the magnitudes.
However, none of these previous approaches started from a notion of `operating condition', related to an {\em asymmetric loss function}. Also, the notion of threshold was not  replaced by a similar concept playing its role for adjusting to the operating condition, and the dual positive-negative character in ROC analysis was blurred.

In this paper we present a graphical representation of regression performance based on a very usual view of operating condition, in regression. Many regression applications have deployment contexts where over-estimations are not equally costly as under-estimations (or vice versa). This is called the {\em loss asymmetry}. 
Loss asymmetry is just a kind of operating condition (or one of its constituents), but a very important one in many applications.

The ROC space for regression is then defined by placing the total over-estimation on the \xaxis and the total under-estimation on the \yaxis. This duality leads to regions and isometrics in the ROC space where over-estimations have less cost than under-estimations and vice versa, and we can plot different regression models to see the notions of dominance. We also consider the construction of hybrid regressors.
The plot leads to curves when we use the notion of {\em shift}, which is just a constant that we can add (or subtract) to example predictions in order to adjust the model to an operating condition. This notion is parallel to the notion of threshold in classification. Interestingly, while we can derive the best shift for a dataset given an existing model (which boils down to shift it to make its average error equal to zero), there are some effective methods to determine this shift for the deployment data given an operating condition, as has been recently explored by \cite{Bansal2008}\cite{zhao2011extended}. Also, there are some other ways to make this shift dependent to each example \cite{probreg2012}. All this leads to a more meaningful interpretation of what the ROC curves in regression mean, and what their areas represent. This will also be explored in this paper.

The paper is organised as follows. 
Section \ref{sec:background} introduces some notation, the problem of context-sensitive evaluation and the use of asymmetric costs in regression. The \RROC space is introduced in section \ref{sec:rocspace}, where we represent several regression models as points, derive the isometrics of the space and develop the notions of hybrid models, dominance and convex hull. Section \ref{sec:roccurves} introduces \RROC curves, which are drawn by ranging a constant shift over the predictions. We introduce an algorithm for plotting them and determine some of its properties in terms of segment slopes and convexity. The area over the \RROC curve (\aoc) is also introduced and analysed.
Section \ref{sec:nonconstant} discusses \RROC curves with non-constant shifts and soft regression models, and the relation with cost plots. Finally, section \ref{sec:conclusions} closes the paper with an enumeration of issues for future investigation.

\section{Context-sensitive problems}\label{sec:background}

In this section we introduce some notation and the basic concepts about context-sensitive regression and the need of asymmetric loss functions.

\subsection{Notation}

Let us consider a multivariate input domain $\mathbb{X} \subset \mathbb{R}^d$ and a univariate output domain $\mathbb{Y} \subset \mathbb{R}$. The domain space $\mathbb{D}$ is then $\mathbb{X} \times \mathbb{Y}$. 
The length of the dataset will usually be denoted by $n$. Examples or instances are just pairs $\left\langle x,y \right\rangle \in \mathbb{D}$, and datasets are subsets of $\mathbb{D}$. 
A {\em crisp} regression model $m$ is a function $m: \mathbb{X} \rightarrow \mathbb{Y}$. A {\em soft regression model} accompanies each prediction with a reliability, confidence or, more generally, a conditional probability density function $\hat{f}(y|x)$ with $y \in  \mathbb{Y}$ and $x \in \mathbb{X}$. 
When the regression model is crisp, we just represent the true value by $y$ and the estimated value by $\hat{y}$. Subindices will be used when referring to more than one example in a dataset.
 
Vectors (unidimensional arrays) are denoted in boldface and its elements with subindices, e.g., $\vect{v}= (v_1,v_2, \dots, v_n)$. Operations mixing arrays and scalar values will be allowed, specially in algorithms, as usual in the matrix arithmetic of many statistical computing languages. For instance, $\vect{v} + c$ means that the constant $c$ is added to all the elements in the vector $\vect{v}$.
The mean of a vector is denoted by $\mu(\vect{v})$ and its standard deviation as $\sigma(\vect{v})$ ---over the population, i.e., divided by $n$.
Given a dataset with $n$ instances $i = 1 \dots n$, the error vector $\vect{e}$ is defined as $e_i \triangleq \hat{y_i} - y_i$. The value  $\mu(\vect{e}^2)$ is known as the mean squared error (\mse), $\mu(\vect{e})$ is known as the mean error (or error bias), $\mu(|\vect{e}|)$ is known as the mean absolute error (\mae) and $\mu(\vect{e})^2$ as the error variance.

\subsection{Context-sensitive problems and loss functions}

In context-sensitive learning \cite{Elk01}, there are several features which describe a context, such as the data distribution, the costs of using some input variables and the loss of the errors over the output variables \cite{turney2000types}. In this paper, we focus on loss functions over the output, which is the kind of costs which ROC analysis deals with (typically integrated, along with the class distribution, within the notion of skew). 
A loss function is defined as follows:

\begin{definition}\label{def:loss}
A loss function is any function 
${\ell}:{\mathbb{Y}} \times {\mathbb{Y}} \rightarrow \mathbb{R}$ which compares elements in the output domain. For convenience, the first argument will be the estimated value, and the second argument the actual value, so its application is usually denoted by ${\ell}(\hat{y},y)$.
\end{definition}

Typical examples of loss functions are the absolute error ($\absloss$) and the squared error ($\squloss$), with $\absloss(\hat{y}, y) = |\hat{y} - y|$ and $\squloss(\hat{y}, y) = (\hat{y} - y)^2$. These two loss functions are {\em symmetric}, i.e. for every $y$ and $r$ we have that ${\ell}(y+r,y)$$=$${\ell}(y-r,y)$. Two of the most common metrics for evaluating regression, the mean absolute error (\mae) and the mean squared error (\mse) are derived from these losses.

\subsection{Asymmetric costs}

Actually, although symmetric loss functions (and derived metrics) are common 
for the evaluation of regression models, it is rarely the case that a real problem has a symmetric cost. For instance, the prediction of sales, consumptions, calls, prices, demands, etc., has almost never a symmetric loss. For instance, a retailing company may need to predict how many items will be sold next week for stock (inventory) management purposes, e.g., in order to calculate how many items must be ordered to refill the stock. Depending on the kind of product, it is usually not the same to over-estimate (increasing stocking costs) than under-estimate (an item is exhausted and it cannot be sold or sold with delays). In fact, it is also rare to find applications where even an asymmetric cost is invariable. For instance, depending on the warehouse saturation, the cost (and the asymmetry) may change in a weekly or daily fashion.
We wish to remark here that a specialised model for a fixed given asymmetry is not the solution in many occasions, either. This motivates the adaptation (or reframing) of models, rather than their re-training for each new asymmetric loss. This is at the core of ROC analysis.

There has been an extensive amount of work on regression using asymmetric loss functions. In some cases, the loss function is embedded in the learning algorithm (see, e.g., \cite{Crone2002,jino2010machine}), which is useful if we know the operating condition during training.  However, the adaptation (or reframing) of an existing model to a different operating condition has also been investigated for regression (e.g., Granger \cite{granger1969prediction,granger1999outline}. 
Many different kinds of asymmetric functions have been explored: {\em Lin-Lin} (asymmetric linear), {\em Quad-Quad} (asymmetric quadratic), {\em Lin-Exp} (approximately linear on one side and exponential on the other side) and {\em Quad-Exp} (approximately quadratic on one side and exponential on the other side) \cite{zellner1986bayesian,christoffersen1996further,christoffersen1997,basu1992bayesian,thompson1996asymmetric}. Some of these approaches try to adapt to the operating condition using complex (generally non-parametric) density functions, which is problematic in general. There are many other approaches. We just mention some of these approaches as an illustration of how important it is in practice to adjust regression models to work with a specific loss function.

As mentioned above, there are many possible asymmetric loss functions. The simplest (and perhaps most common) one is the asymmetric absolute error $\aabsloss$:

\begin{definition} \label{def:absloss}
The asymmetric absolute error $\aabsloss$ is a loss function defined as follows:
\begin{eqnarray*}
 \aabsloss(\hat{y}, y) \:\: = & 2\alpha (y - \hat{y}) & \mbox{if} \:\:\: \hat{y} < y \\
                                   = &  2(1-\alpha)(\hat{y} - y) & \mbox{otherwise}
\end{eqnarray*}						
\end{definition}

with $\alpha$ being the cost proportion (or asymmetry) between 0 and 1, with increasing values meaning higher cost for {\em low predictions} (underestimation). In other words, when $\alpha=0$ we mean that predictions below the actual value have no cost. When $\alpha=1$ we mean that predictions above the actual value have no cost.
When $\alpha=0.5$ we mean that costs above and below are symmetric.


\section{The \RROC space}\label{sec:rocspace}

For every regression model deployed to a new dataset we can determine the error for each example and whether it is an over-estimation or under-estimation. 
More formally:

\begin{definition}The total over-estimation is given by $\OVER \triangleq \sum_i \{e_i  \:\:|\:\:   e_i > 0\}$ and the total under-estimation is given by $\UNDER \triangleq \sum_i \{e_i  \:\:|\:\:   e_i < 0\}$. 
\end{definition}

The following example illustrate this:

\begin{example}\label{ex:example1}
Consider a regression model $m_1$ which is applied to a dataset with $n=10$ examples $e_1 \dots e_{10}$, issuing the predicted values $\hat{y}$ and actual values $y$:

\vspace{0.4cm}
{	\centering
{\footnotesize
\noindent		\begin{tabular}{r|cccccccccc}
		& $1$ & $2$ & $3$ & $4$ & $5$ & $6$ & $7$ & $8$ & $9$ & ${10}$ \\
	$\hat{y}$ & -0.082 & 3.323 & 2.320 & 1.080 & 7.893 & 4.983 & 5.121 & 3.442 & 2.083 & 1.112	\\	
	$y$ & 0.211 & 2.725 & 1.933 & 3.242 & 7.858 & 6.061 & 7.173 & 3.082 & 0.894 & 1.203 \\
	$e$ &  -0.293 & 0.598 & 0.387 & -2.162 & 0.035 & -1.078 & -2.052 & 0.360 & 1.189 &-0.091 \\
		\end{tabular}
}
}
\vspace{0.4cm}

The error row ($e$) shows the difference, which is positive for over-estimations and negative for under-estimations. The sum of over-estimations ($\OVER$) is $2.569$ while the sum of under-estimations ($\UNDER$) is $-5.676$. This regression model clearly under-estimates (it has a negative error bias, since $\mu(\vect{e}) < 0$). The \mae ($0.825$) and the \mse ($1.219$) do not show the asymmetry of predictions.
\end{example}

\subsection{Showing models in \RROC space}

Certainly, different regression models would show different error asymmetries (or error bias). The basic idea of the ROC space for regression is to show this asymmetry:

\begin{definition}
The Regression Receiver Operating Characteristic (\RROC) space is defined as a plot where we depict total over-estimation ($\OVER$) on the \xaxis and total under-estimation ($\UNDER$) on the \yaxis. Since $\OVER$ is always positive (but unbounded) and $\UNDER$ is always negative (but unbounded), we typically will place the point $(0,0)$ on the upper left corner (the \RROC heaven), and will clip both the \xaxis and \yaxis as necessary to show the region of interest.
\end{definition}

\begin{figure}
\centering
\includegraphics[width=0.5\textwidth]{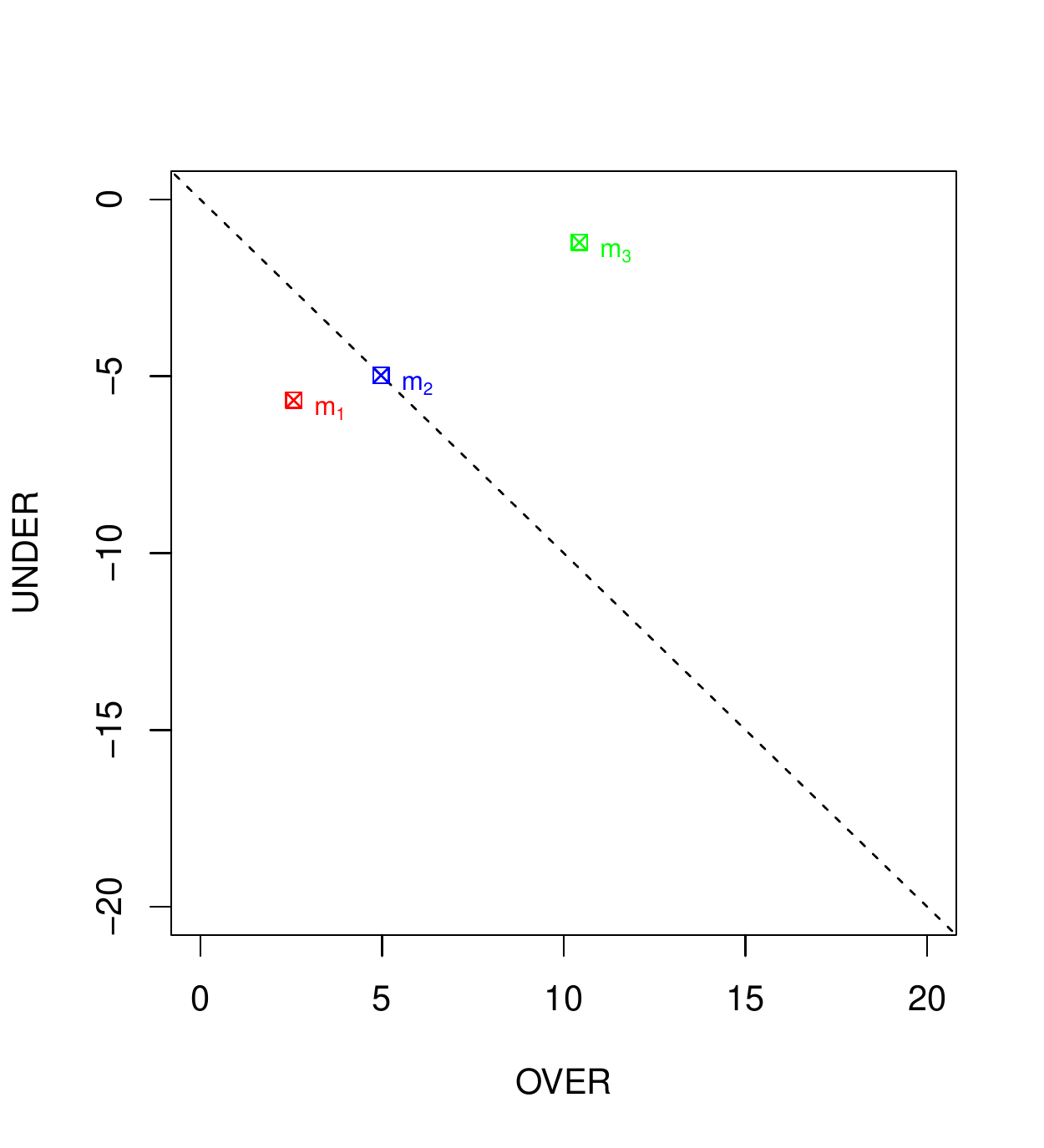} 
\vspace{-0.5cm}
\caption{\RROC space and the representation for regression models $m_1$ (in red), $m_2$ (in blue) and $m_3$ (in green) in examples \ref{ex:example1}, \ref{ex:example2} and \ref{ex:example3}. The diagonal (dashed) shows where \UNDER and \OVER are equal. Model $m_2$ has zero error bias ($\mu(\vect{e}) = 0$). }
\label{fig:RROCnoshift}
\end{figure}

Figure \ref{fig:RROCnoshift} shows the \RROC space and the regression model $m_1$ in example \ref{ex:example1}. We will occasionally draw a diagonal line $\OVER-\UNDER=0$ to show the points where the under-estimation equals the over-estimation.

One might argue why we use absolute values for the \xaxis and \yaxis instead of relative values. In fact, ROC analysis uses relative values. There are two reasons for this. First, using relative values would not make the \RROC space finite. Second, and more importantly, using relative values we could have cases where changing a single infinitesimal change on one example could end up at very different locations. For instance, consider the error vectors $e_A=\{-10, -0.1, 5\}$ and $e_B=\{-10, 0.1, 5\}$. While \UNDER and \OVER are almost the same, the relative \UNDER and \OVER would be $\{-5.05, 5\}$ and $\{-10, 2.55\}$ for two almost equal error vectors. 
This justifies that the \RROC space shows absolute values. In this sense, and strictly speaking, the parallel with ROC analysis for classification can be done with the `coverage curves' \cite{petersbook}, which are the absolute variant of ROC curves.

Let us now consider a second model:

\begin{example}\label{ex:example2}
Consider a regression model $m_2$ which is applied to the same dataset as example \ref{ex:example1}:

\vspace{0.4cm}
{	\centering
{\footnotesize
\noindent		\begin{tabular}{r|cccccccccc}
		& $1$ & $2$ & $3$ & $4$ & $5$ & $6$ & $7$ & $8$ & $9$ & ${10}$ \\
	$\hat{y}$ & 0.786 & 2.078 & 0.587 & 1.676 & 9.052 & 5.875 & 6.885 & 3.038 & 4.097 & 0.308 \\	
	$y$ & 0.211 & 2.725 & 1.933 & 3.242 & 7.858 & 6.061 & 7.173 & 3.082 & 0.894 & 1.203 \\
	$e$ &  0.575 & -0.647 & -1.346 & -1.566 & 1.194 & -0.186 & -0.288 & -0.044 &  3.203 & -0.895\\
		\end{tabular}
}
}
\vspace{0.4cm}

The sum of over-estimations ($\OVER$) is $4.972$ while the sum of under-estimations ($\UNDER$) is $-4.972$. This regression finds an equilibrium between over and under-estimations (it is unbiased, since $\mu(\vect{e}) = 0$). The \mae ($0.9944$) and the \mse ($1.7619$) are worse than $m_1$ in example \ref{ex:example1}.
\end{example}

This model ($m_2$) with  $\OVER-\UNDER=0$ is also shown in Figure \ref{fig:RROCnoshift}. Clearly it is on the diagonal.

Finally let us consider a third model:

\begin{example}\label{ex:example3}
Consider a regression model $m_3$ as follows:

\vspace{0.4cm}
{	\centering
{\footnotesize
\noindent		\begin{tabular}{r|cccccccccc}
		& $1$ & $2$ & $3$ & $4$ & $5$ & $6$ & $7$ & $8$ & $9$ & ${10}$ \\
	$\hat{y}$ & 1.253 & 4.232 & 1.734 & 5.325 & 6.842 & 9.325 & 8.232 & 3.525 & 1.352 & 1.778 \\
	$y$ & 0.211 & 2.725 & 1.933 & 3.242 & 7.858 & 6.061 & 7.173 & 3.082 & 0.894 & 1.203 \\
	$e$ &   1.042 &   1.507 &  -0.199 &   2.083 &  -1.016 &   3.264 &   1.059 &   0.443 &   0.458 &   0.575 \\
		\end{tabular}
}
}
\vspace{0.4cm}

In this case, the sum of over-estimations ($\OVER$) is $10.431$ while the sum of under-estimations ($\UNDER$) is $-1.215$. This regression model clearly over-estimates (it has a positive error bias, since $\mu(\vect{e}) > 0$). The \mae ($1.165$) and the \mse ($2.12$) show that this model is, in terms of overall error, worse than models $m_1$ and $m_2$.
\end{example}

From each point in \RROC space, we can derive its \mae very easily. For model $m_3$, for example, we have that $\mae=1.165 = (\OVER - \UNDER) / n$, so it is just half the perimeter of the rectangle that each point creates with the \RROC heaven (0,0). In other words, the \mae (more precisely the absolute error) is just the Manhattan distance to \RROC heaven.
It is important to note that the diagonal (the Euclidean distance) is just given by $\sqrt(\OVER^2 + \UNDER^2)$, which we call $\mmse$ (as a macro-averaged version of $\mse$). This $\mmse$ measure is interesting in itself, because highly penalises models for which there is a high imbalance in over and under-estimations, and can be seen, in some way, as a measure of `symmetric calibration' \arxiv{\cite{Handbook}}.

In \RROC space we denote the regression model always outputting $\infty$ and the model always outputting $-\infty$ as the (trivial) {\em extreme} regression models, which fall at $(\infty, 0)$ and $(0, -\infty)$ respectively in \RROC space.

\subsection{\RROC space isometrics}

We have mentioned above that ($1/2$ of) the perimeter of the rectangle from \RROC heaven to the regression model corresponds to $\mae$. Can we extend this observation to the asymmetric loss?
The following straightforward lemma shows that total asymmetric absolute loss can be calculated graphically as the sum of the distance to the \yaxis ($\OVER = 0$) and to the \xaxis ($\UNDER=0$), using the appropriate asymmetry factor $\alpha$.

\begin{lemma}\label{lem:perimeter}
The total asymmetric absolute loss is given by:
\[L = \sum_i\aabsloss(\hat{y_i}, y_i) = - 2\alpha \cdot \UNDER + 2(1-\alpha) \cdot \OVER. \]
\end{lemma}
\begin{proof}
\begin{eqnarray*}
L = \sum_i\aabsloss(\hat{y_i}, y_i) & = & \sum_i \{ 2\alpha (y_i - \hat{y_i})  \mbox{\:\:\:if\:\:\:} \hat{y_i} < y_i, 2(1-\alpha)(\hat{y_i} - y_i) \mbox{\:\:\:otherwise} \} \\
& = &  \sum_i \{2\alpha(-e_i)  \:\:|\:\:   e_i < 0 \} + \sum_i \{ 2(1-\alpha)(e_i)  \:\:|\:\:   e_i > 0 \} \\ 
& = &  -2\alpha \cdot \UNDER + 2(1-\alpha) \cdot \OVER
\end{eqnarray*}
\end{proof}

Clearly, for $\alpha=0.5$, we have that this is the absolute error.
All this also shows that the closer we are to \RROC heaven $(0,0)$ (in terms of a Manhattan distance) the better. Finally, this leads to loss isometrics:

\begin{definition}\RROC isometrics are defined by varying $t$ over:
\[ -2\alpha \cdot \UNDER + 2(1-\alpha) \cdot \OVER = t \]
\end{definition}

We can get any of the infinite (and parallel) isometrics. The following proposition just gets the slope of each isometric:

\begin{proposition}\label{prop:slope}
Given an isometric $-2\alpha \cdot \UNDER + 2(1-\alpha) \cdot \OVER = t$, the slope only depends on $\alpha$ and is given by:
\[ slope = \frac{1-\alpha}{\alpha} \]
\end{proposition}
\begin{proof}
By isolating the variable $\UNDER$ we have:
\[ \UNDER = \frac{t -  2(1-\alpha) \cdot \OVER}{-2\alpha} = \frac{-2t}{\alpha} + \frac{1-\alpha}{\alpha} \cdot \OVER\]
The $slope$ is then given by the second term $\frac{1-\alpha}{\alpha}$
\end{proof}

Clearly, for $\alpha=0$ (under-estimations have no cost) and we have infinite slope. For $\alpha=1$ (over-estimations have no cost), we would have a slope 0.

This notion of isometric is very similar to the notion already present in ROC analysis for classification \cite{Fla03}. In fact, this means that we can slide isometrics to find optimal points in \RROC space, in the very same way as we do in ROC space.

\begin{figure}
\centering
\includegraphics[width=0.5\textwidth]{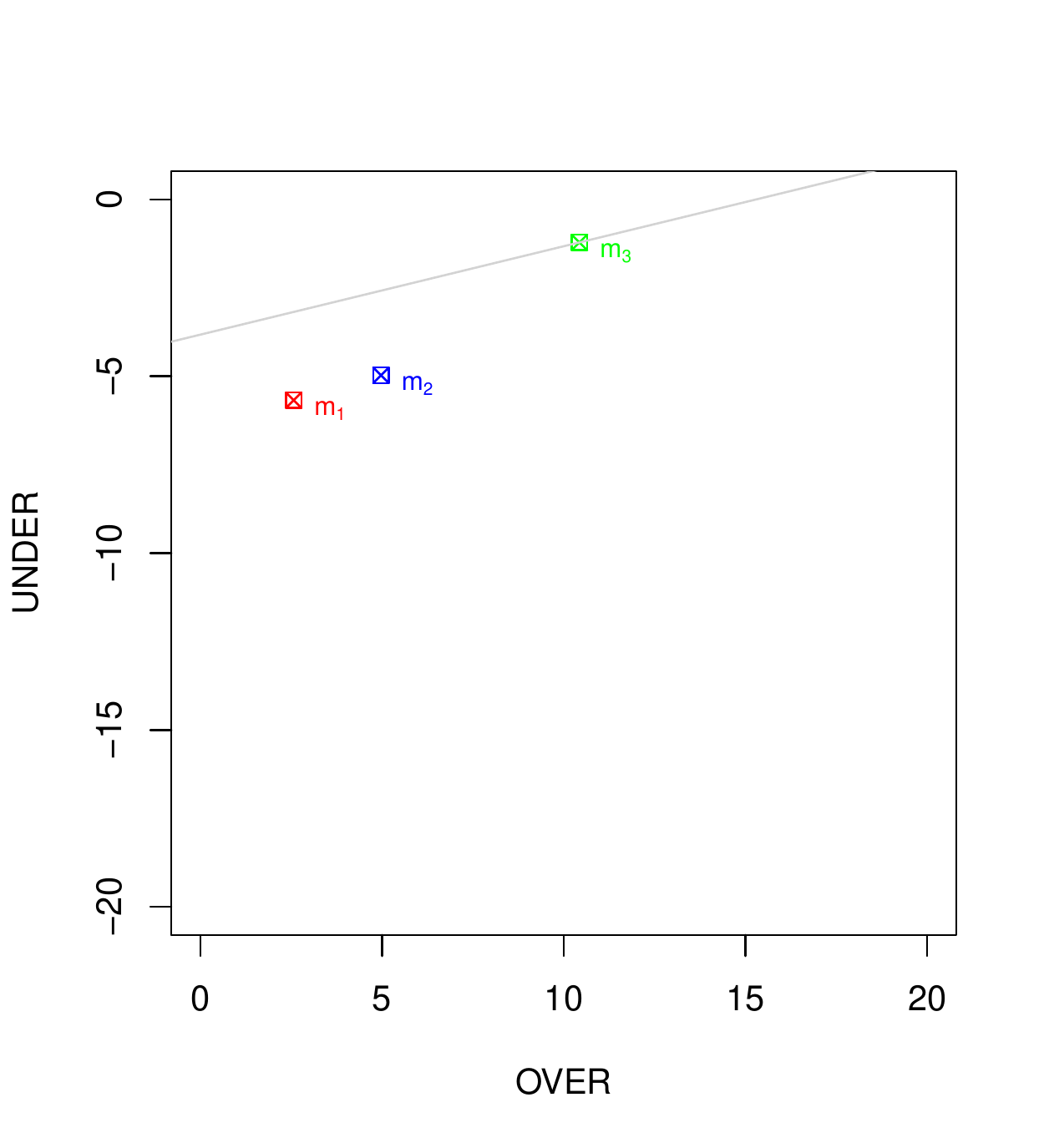} 
\vspace{-0.5cm}
\caption{The three models as in Figure \ref{fig:RROCnoshift}. We show the first isometric line (light grey) corresponding to $\alpha= 0.8$ ($slope = 0.25$) touching any of the three models. }
\label{fig:RROCnoshiftiso}
\end{figure}

Let us illustrate this.
Figure \ref{fig:RROCnoshiftiso} shows the \RROC space and the regression models $m_1$, $m_2$ and $m_3$ in examples \ref{ex:example1}, \ref{ex:example2}  and \ref{ex:example3} respectively. We also consider the operating condition $\alpha = 0.8$, meaning that under-estimations are 4 times more expensive than over-estimations. This $\alpha$ leads to a slope of $0.25$. By sliding through all the parallel isometric lines from the one crossing the \RROC heaven $(0,0)$ to the first isometric touching a point corresponding to any model, we touch at $(10.431,-1.215)$ first.
In fact, the $intercept$ is given by isolating it from the line equation $under = slope \cdot \OVER + intercept$, i.e., $intercept = \UNDER - slope \cdot \OVER$, which, in this case, leads to $-3.82275$. The line $\UNDER = 0.25 \cdot \OVER -3.82275$ is then shown on Figure \ref{fig:RROCnoshiftiso}, touching regression model $m_3$.
Even though model $m_3$ has a worse mean (symmetric) absolute error than $m_1$, for this operating condition $\alpha$, it leads to lower total asymmetric absolute error. While $m_1$ has a loss of $-2\alpha \cdot \UNDER + 2(1-\alpha) \cdot \OVER = -1.6 \cdot (-5.676) + 0.4 \cdot 2.569 = 10.1092$, we have that $m_3$ has a loss of $-2\alpha \cdot \UNDER + 2(1-\alpha) \cdot \OVER = -1.6 \cdot (-1.215) + 0.4 \cdot 10.431 = 6.1164$.

\subsection{Hybrid models, dominance and convex hull}

Another construction that is also originally present in ROC analysis for regression is the notion of hybrid models.
Given any two models, we can construct a hybrid model by randomly choosing each prediction from any of both models using a (biased) coin. Note that this is very different to averaging both models.

\begin{figure}
\centering
\includegraphics[width=0.5\textwidth]{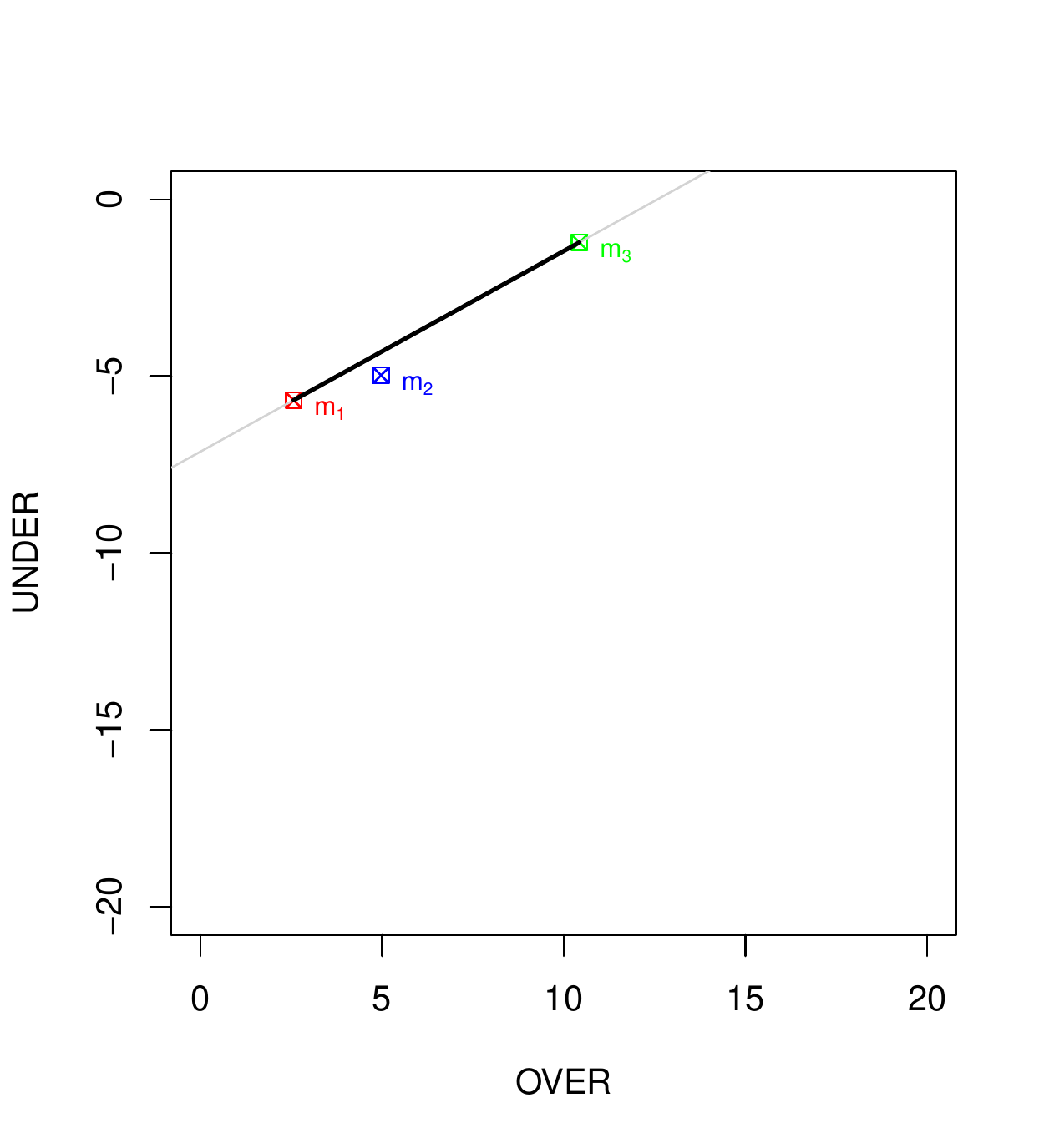} 
\vspace{-0.5cm}
\caption{The three models as in Figure \ref{fig:RROCnoshift}. By considering any model which can be constructed by just choosing predictions randomly (with any bias) between models $m_1$ and $m_3$, we can show a segment of models (in solid black). }
\label{fig:RROCnoshift13}
\end{figure}

Figure \ref{fig:RROCnoshift13} shows the isometric (in light grey) passing through models $m_1$ and $m_3$. The solid black segment connecting both models shows that any model along the segment can be constructed. More precisely, each point in that segment would represent the expected value of a model constructed in this way. Consequently, we can just connect both points since any point in between is technically achievable (at least in expectation). 

In this particular case, we just draw a line between the point representing $m_1$: $(2.569,-5.676)$ and the point representing $m_3$: $(10.431,-1.215)$, leading to 
$\UNDER = 0.567 \cdot \OVER - 7.134$.
%
%
From this slope of $0.567$, we just calculate 
%
$\alpha = \frac{1}{1+slope} = 0.638$.
Obviously, for this $\alpha$ both models have the same loss. $L(m_1) = 0.638 \cdot 5.676 + (1-0.638) \cdot 2.569 = 4.551$ and $L(m_3) =0.638 \cdot 1.215 + (1-0.638) \cdot 10.431 = 4.551$.

Given these two models, we say that, for slopes lower than $0.567$ and asymmetries $\alpha$ greater than $0.638$, model $m_3$ dominates, while we have that model $m_1$ dominates for the rest of operating conditions.

\begin{figure}
\centering
\includegraphics[width=0.5\textwidth]{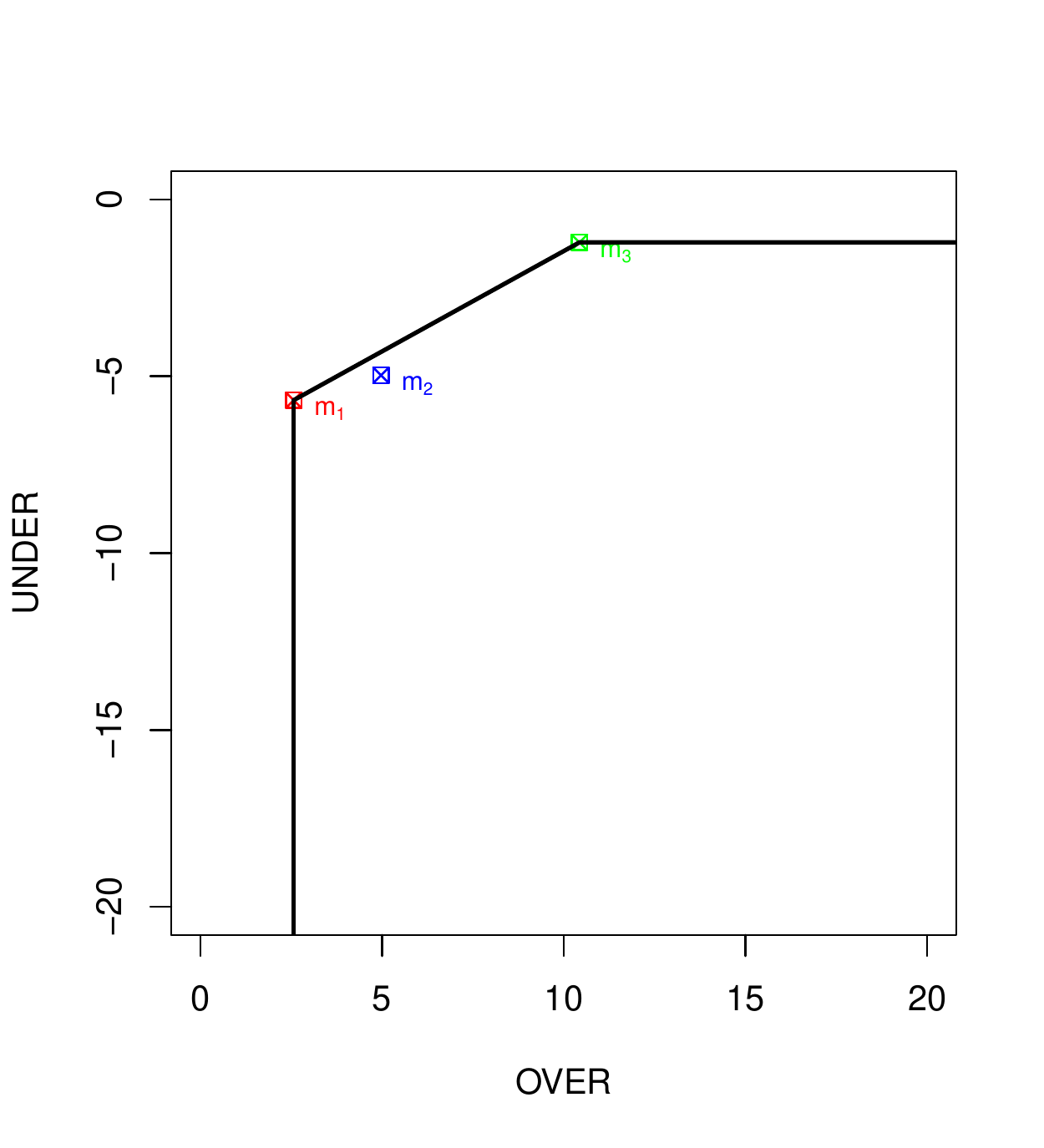} 
\vspace{-0.5cm}
\caption{The three models as in Figure \ref{fig:RROCnoshift}. By considering any model which can be constructed by just choosing predictions randomly (biassedly) between any other two models ---including the extreme models at $(0,-\infty)$ and $(\infty, 0)$--- we can derive the convex hull (shown in solid black). }
\label{fig:RROCnoshifthull}
\end{figure}

This leads to the notion of dominance and convex hull. In fact, when connecting all the points by the segments representing the hybrid models (and also including the extreme classifiers at $(0,-\infty)$ and $(\infty, 0)$, we can calculate the convex hull, since any model under the convex hull can be discarded, in the same way as traditional ROC analysis DOES.
Figure \ref{fig:RROCnoshifthull} shows the convex hull of the three models and the extreme models.
We see that model $m_2$ can be discarded. It cannot be optimal for any operating condition.

\section{\RROC curves}\label{sec:roccurves}

In ROC analysis for classification, we can tweak the predictions of a crisp classifier by changing the predicted class to a random percentage of examples. With this, we can move the classifier in the ROC space, but this just moves the classifier along the two straight lines that connect the original point with the points at $(0,0)$ and $(1,1)$ (the trivial, or extreme, classifiers). For this reason, occasionally a crisp classifier is represented in ROC space as trapezium\footnote{A trapezoid in American English.}, connecting the point which corresponds to the classifier with the extreme classifiers. This two-segment `curve' does not bring more information than the original point, but shows that other TPR and FPR can be achieved by this random swapping of examples. In the end, it just shows the hybrid classifier constructed with the extreme classifiers.

In general, however, in ROC analysis, curves are constructed by the use of soft classifiers, i.e., classifiers which output a rank, score or probability estimation. By moving a threshold from the lowest possible valuable to the highest possible value (or vice versa) we get many possible crisp classifiers, each of them represented by a point in ROC space.

Interestingly, in \RROC space, we do not need soft regression models in order to create a {\em curve}. It is just sufficient to use a {\em shift}, which works as a parallel concept to the notion of threshold. For each example we can get a modified prediction as $\dot{y} \leftarrow \hat{y} + s$, where $s$ is the shift. Although there are, as we will see, many ways of determining this shift, it seems natural to consider first that $s$ is constant, i.e., that we apply the same value for all the examples. 

\begin{definition}
Given a regression model $m$, a (constant-)shifted regression model, denoted by $m\left\langle s \right\rangle$, is the result of adding the same shift $s$ to all its predictions, i.e., $\hat{y}' \leftarrow \hat{y} + s$ for all predictions $\hat{y}$.
\end{definition}

This shift can be moved from the lowest possible value ($-\infty$) to the maximum possible value ($\infty$). 
%
%
 This leads to the notion of \RROC curve.

\begin{definition}
Given a regression model $m$, its \RROC curve using a (constant) shift is given by plotting all the models $m\left\langle s \right\rangle$ with $s$ ranging in $[-\infty,\infty]$.
\end{definition}

We can instantly plot the curves pointwise, by just using a sufficient dense range of values for $s$. However, there is a more direct way of plotting and analysing the \RROC curve if we investigate a little bit. This is what we do next.

\subsection{Algorithm for drawing \RROC curves}

We can realise that if we move the shift from $s_1$ to $s_2$ and no example changes from \OVER to \UNDER or vice versa, then the increment/decrement in \OVER and \UNDER is linear, as the following proposition shows: 

\begin{proposition}\label{prop:line}
Given a model $m$, for any two shifts $s_1$ and $s_2$ such that the examples for which $m\left\langle s_1 \right\rangle$ and $m\left\langle s_2 \right\rangle$ over-estimate are the same (and hence the rest that under-estimate are also the same for both), then for any other shift $s_3$ with $s_1 \leq s_3 \leq s_2$ we have that the points $(\OVER,\UNDER)$ for the three models  $m\left\langle s_1 \right\rangle$, $m\left\langle s_2\right\rangle$ and $m\left\langle s_3 \right\rangle$ lie on the same straight line.
\end{proposition}

\begin{proof}
We have that $\OVER$ for $m\left\langle s_1 \right\rangle$ is calculated as:
$\OVER_1 = \sum_i \{e_i + s_1  \:\:|\:\:   e_i + s_1 > 0 \}$ while  $\OVER$ for $m\left\langle s_2 \right\rangle$ is calculated as:
$\OVER_2 = \sum_i \{e_i + s_2  \:\:|\:\:   e_i + s_2 > 0 \}$. 
Since, by assumption, the examples which over-estimate are the same for $m\left\langle s_1 \right\rangle$ and $m\left\langle s_2 \right\rangle$, let us call this number $n_o$. The previous two expressions can then be rewritten as:
\[ \OVER_1 = n_o s_1 + \sum_i \{e_i  \:\:|\:\:   e_i + s_1 > 0\} \]
\[ \OVER_2 = n_o s_2+ \sum_i \{e_i  \:\:|\:\:   e_i + s_1 > 0\} \]
Note that the second term is also rewritten with $s_1$, since the elements are the same. 
In this way, we express that the second term is equal. 
Also, since the examples which over-estimate are the same for $s_1$ and $s_2$ they have to be the same necessarily for every $s_3$ with $s_1 \leq s_3 \leq s_2$ as well. So, we also have:
\[ \OVER_3 \triangleq n_o s_3+ \sum_i \{e_i  \:\:|\:\:   e_i + s_1 > 0\}\]
We can see that these three co-ordinates only differ on the first term, which is linearly related to $s$ ($s_1$, $s_2$ or $s_3$). 
We can obtain similar expressions for $\UNDER_1$, $\UNDER_2$ and $\UNDER_3$ and their $n_u$ examples.
This means that the three points are related by a linear term on $s$, expressed as $(n_o s, n_u s)$ so they lie on the same line.
\end{proof}

From proposition \ref{prop:line} we can introduce a very simple algorithm to draw \RROC curves:

\IncMargin{1em}
\begin{algorithm}
\DontPrintSemicolon
\SetKwInOut{Input}{input}\SetKwInOut{Output}{output}
\TitleOfAlgo{PlotRROCCurve} 
\Input{Two arrays $\vect{\hat{y}}$ and $\vect{y}$ of size $n$ with the predicted and the actual values respectively.}
\Output{The $n+2$ vertex points of the \RROC Curve in arrays $\vect{RROCX}$ and $\vect{RROCY}$}
\BlankLine
{// Draws the curve from bottom-left corner to top-right corner}\;  
  $\vect{e} \leftarrow  \vect{\hat{y}} - \vect{y}$ \;
	$\vect{e} \leftarrow SortDecreasingly(\vect{e})$ \; 
  $RROCX_1 \leftarrow 0$ \;
  $RROCY_1 \leftarrow -\infty$ \;
  \For{$i\leftarrow 1$ \KwTo $n$}{
	  $s \leftarrow e_i \:\:\:\:$ { // The shift $s$ as examples change from \OVER to \UNDER} \; 
	  $\vect{t} \leftarrow \vect{e} - s \:\:\:\:$   { // Applies a constant shift $s$ to the array $\vect{e}$} \; 
    $RROCX_{i+1} \leftarrow \sum_j \{ t_j \:\:  |  \:\: t_j > 0    \}    \:\:\:\:$ // \OVER  \;
    $RROCY_{i+1} \leftarrow \sum_j \{ t_j \:\: | \:\: t_j \leq 0 \}    \:\:\:\:$ // \UNDER \;
  }
  $RROCX_{n+2} \leftarrow \infty$ \;
  $RROCY_{n+2} \leftarrow 0$ \;
	\BlankLine
	\caption{Algorithm for drawing a \RROC curve. We use brackets for array notation and array operations. The algorithm can be further simplified by updating the array $\vect{e}$ and calculating $\OVER$ and $\UNDER$  incrementally in each iteration in the loop.}\label{algo:RROCcurve}
\end{algorithm}
\DecMargin{1em}

From the first line of the algorithm, we see that the \RROC Curve can be drawn by just giving the error vector (e.g., the last row in examples \ref{ex:example1}, \ref{ex:example2} and \ref{ex:example3}).

\begin{figure}
\centering
\includegraphics[width=0.5\textwidth]{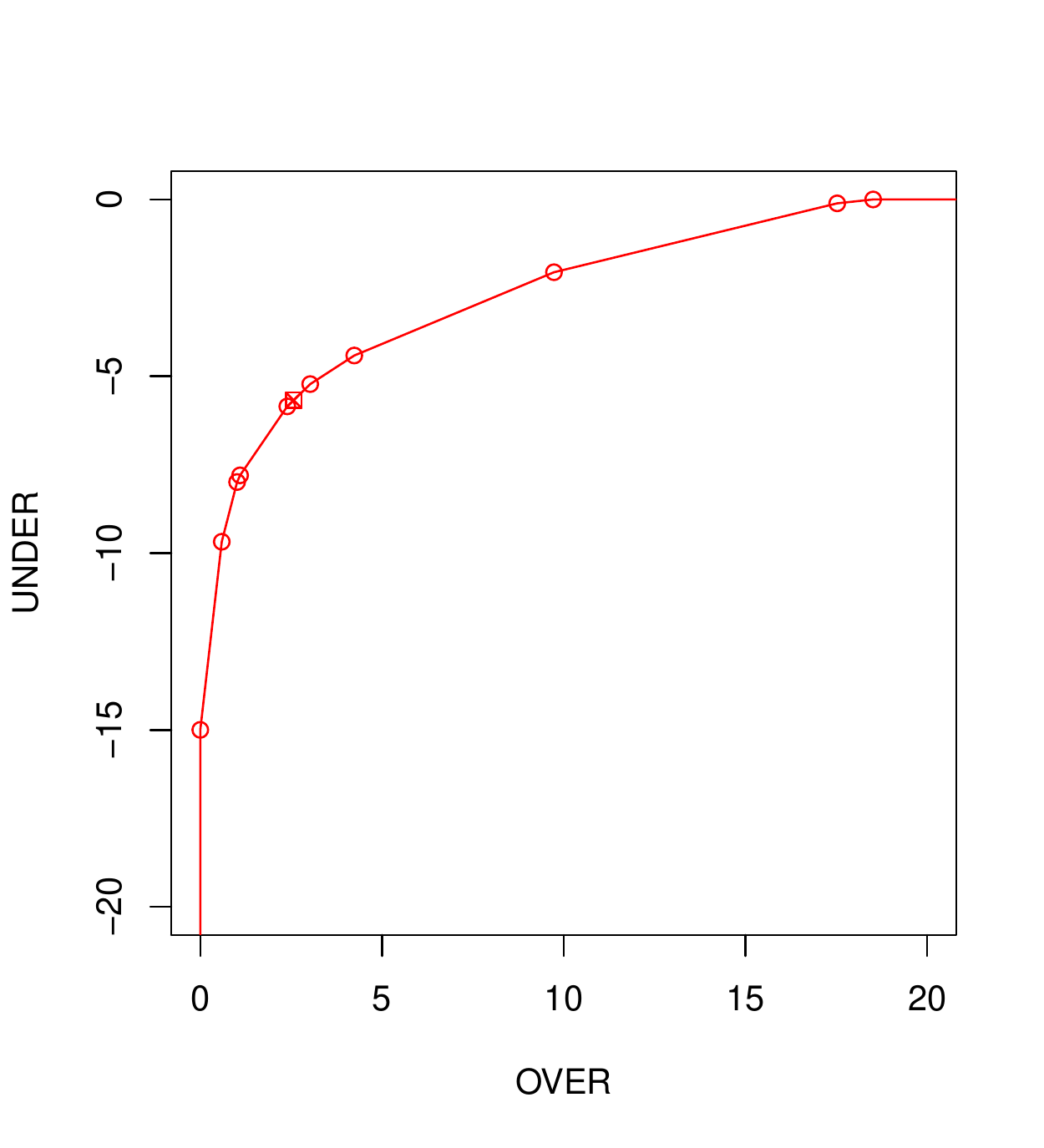} 
\vspace{-0.5cm}
\caption{Model $m_1$ in example \ref{ex:example1} drawn as a \RROC curve by changing the shift. Vertex points ($10$ in this case, since the two extremes are not visible in the plot) are shown as small circles. The curve is then composed of $11$ segments (there are $10$ examples).  The original shift ($s=0$) is still represented with a small square and lies on a segment between two vertex points.}
\label{fig:RROCshift1}
\end{figure}

Figure \ref{fig:RROCshift1} shows a \RROC curve using this algorithm for $m_1$ in example \ref{ex:example1}.
The points where the slope of the \RROC Curve change are called {\em vertex points}, and the rest of points are said to fall onto the segments. Consequently a \RROC Curve for a regression model applied to a dataset with $n$ instances has $n+2$ vertex points (typically, only $n$ are visible on the plot, because two are the extreme points) and $n+1$ segments, denoted by $\overline{i,i+1}$ with $i=1 \dots n+1$. We clearly see $n=10$ points on Figure \ref{fig:RROCshift1}.

\begin{figure}
\centering
\includegraphics[width=0.5\textwidth]{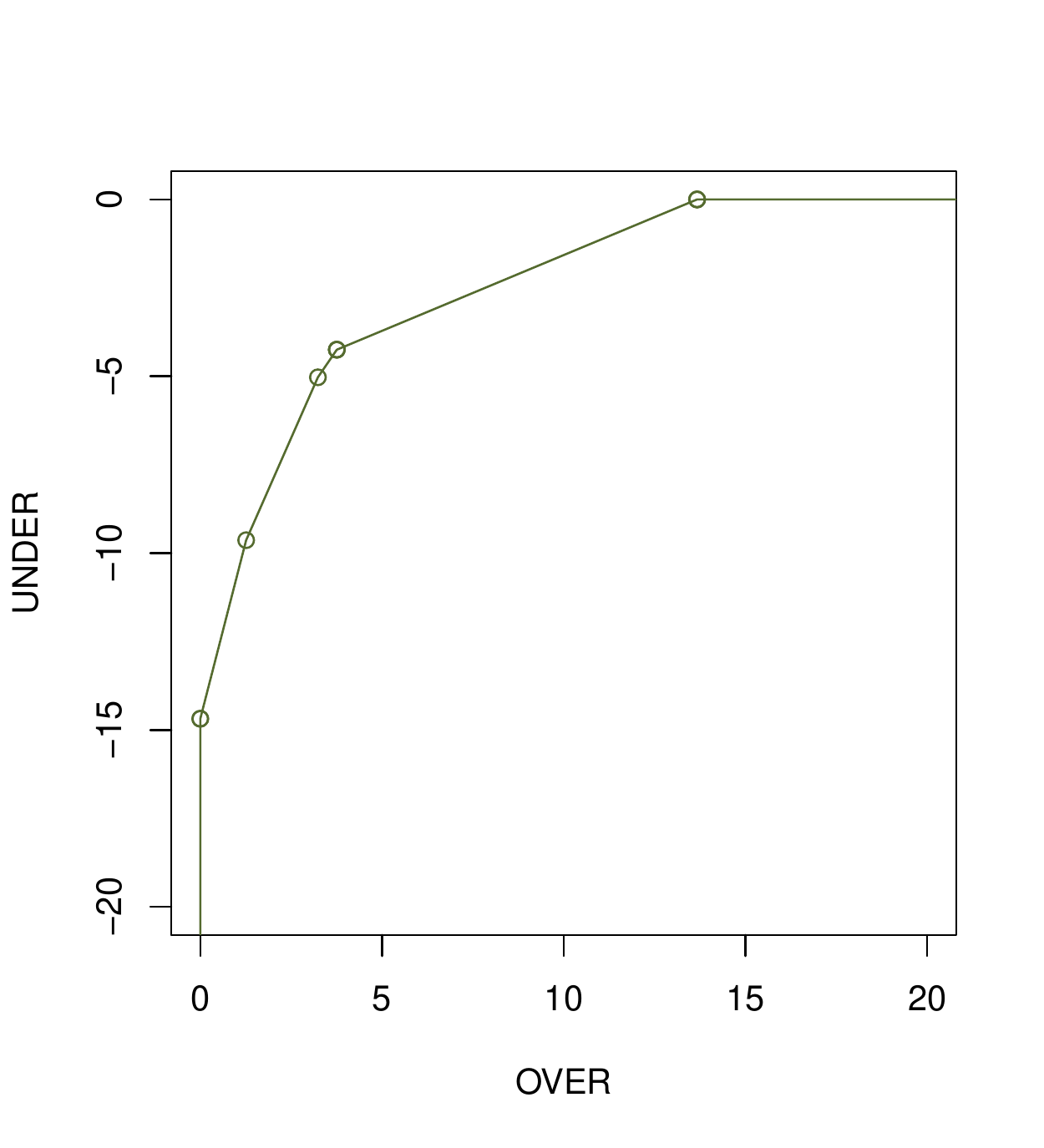} 
\vspace{-0.5cm}
\caption{Model $m_4$ in example \ref{ex:example4} drawn as a \RROC curve by changing the shift. The model has several errors with the same value (two triple ties and a double tie), so the number of distinct visible points is reduced from $n=10$ to $n-2-2-1= 5$.}
\label{fig:RROC4ties}
\end{figure}

In case there are some ties in the error vector, then some of these vertex points and segments collapse into a single point. Figure \ref{fig:RROC4ties} shows

\begin{example}\label{ex:example4}
Consider a regression model $m_4$ as follows:

\vspace{0.4cm}
{	\centering
{\footnotesize
\noindent		\begin{tabular}{r|cccccccccc}
		& $1$ & $2$ & $3$ & $4$ & $5$ & $6$ & $7$ & $8$ & $9$ & ${10}$ \\
	$\hat{y}$ & 0.123 & 1.221 & 1.845 & 4.573 & 8.558 & 7.392 & 5.669 & 1.578 & 0.806 & 1.245	\\	
	$y$ & 0.211 & 2.725 & 1.933 & 3.242 & 7.858 & 6.061 & 7.173 & 3.082 & 0.894 & 1.203 \\
	$e$ &   -0.088 & -1.504 & -0.088 & 1.331 & 0.700 & 1.331 & -1.504 & -1.504 & -0.088 & 0.042 \\
		\end{tabular}
}
}

\vspace{0.4cm}
We see a triple tie between examples 1, 3 and 9, another triple tie between examples 2, 7 and 8, and a double tie between examples 4 and 6. With this, there are only 5 different error values.
\end{example}

\subsection{Properties: slope and convexity}

From the new \RROC curve, we may want to determine the slopes of each segment, in order to exactly determine where each possible isometric (and asymmetry $\alpha$) would lead to on the curve. This can be done very easily, as the following lemma shows:

\begin{lemma}\label{lemma:slope}
The slope of each segment $\overline{i,i+1}$ in the \RROC curve is given by $(n+1-i)/(i-1)$, with $i=1 \dots n+1$.
\end{lemma}
\begin{proof}
Let us assume no ties in the error vector.
As shown in proposition \ref{prop:line}, there is one example changing from \UNDER to \OVER (from bottom-left to top-right) at each vertex point. At the first vertex point $i=1$, all the examples are under-estimated, and the shift change moves along an infinite slope. For the next vertex point $i=2$, we have $n-1$ under-estimated examples and $1$ over-estimated example. This means that the shift change moves along one unit right and $n-1$ units up, with a slope of $n-1$. By induction, this leads to $(n+1-i)/(i-1)$, with the last segment having 0 slope.
If there are ties, the result is similar with more than one example changing from under-estimation to over-estimation at a time.
\end{proof}

Thus, and somewhat surprisingly, given a fixed number of examples, several regression models will have exactly the same slopes. The difference between the curves will be given by the length of the segments, not their slopes.
From the equation $\frac{1-\alpha}{\alpha}= slope$ in proposition \ref{prop:slope} relating asymmetries and slopes, we have that each segment $\overline{i,i+1}$ corresponds to an $\alpha =  \frac{1}{slope+1}$, leading to $\alpha=\frac{i-1}{n}$ with $i=1 \dots n+1$.

Finally, from the previous Figure \ref{fig:RROCshift1}, we see that the curve is convex. Is this true in general? The following proposition shows it is.

\begin{proposition}\label{prop:convex}
For every regression model, the \RROC Curve is convex\footnote{Note that in ROC analysis we typically say `convex' when the region below is a convex set, while, generally, in mathematics, this refers to the region above being a convex set.}.
\end{proposition}
\begin{proof}
It is direct from lemma \ref{lemma:slope} since the sequence of the segment slopes of the curve $(n+1-i)/(i-1)$ is non-increasing.
\end{proof}

\begin{figure}
\centering
\includegraphics[width=0.5\textwidth]{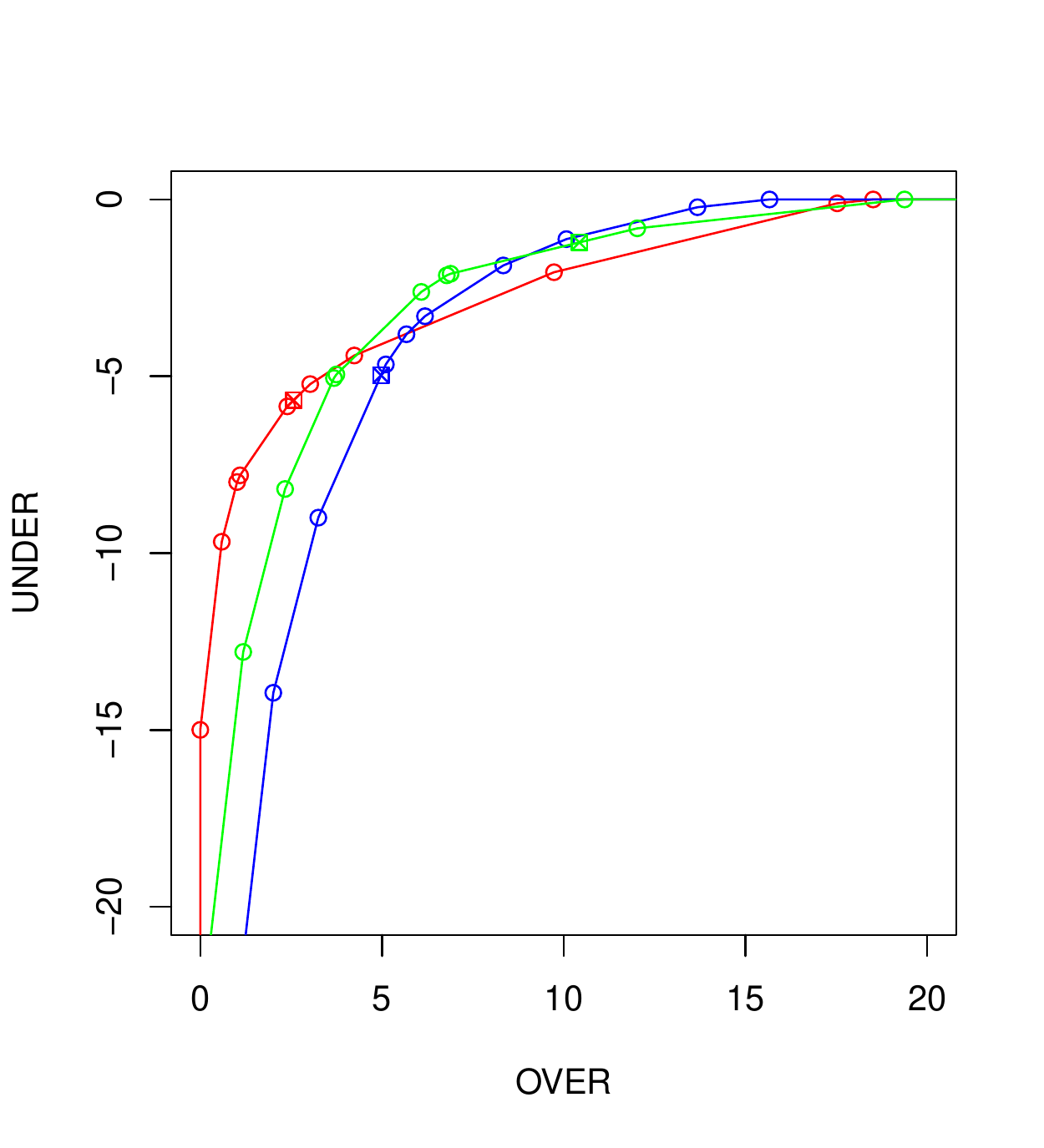} 
\vspace{-0.5cm}
\caption{The three models $m_1$ (red), $m_2$ (blue) and $m_3$ (green) in examples \ref{ex:example1}, \ref{ex:example2} and \ref{ex:example3} drawn as \RROC curves by changing the shift. Note that in this case model $m_3$ cannot be rejected, because there are regions where it is optimal. If we select the best portions from the three models we see concavities, which can be resolved by the use of a convex hull.}
\label{fig:RROCshiftall}
\end{figure}

The convexity of a single \RROC curve does not mean that the notion of convex hull seen in the previous section is useless for curves. More on the contrary. Whenever we have more than one model, we can see concavities. Figure \ref{fig:RROCshiftall} precisely shows this.

\begin{figure}
\centering
\includegraphics[width=0.5\textwidth]{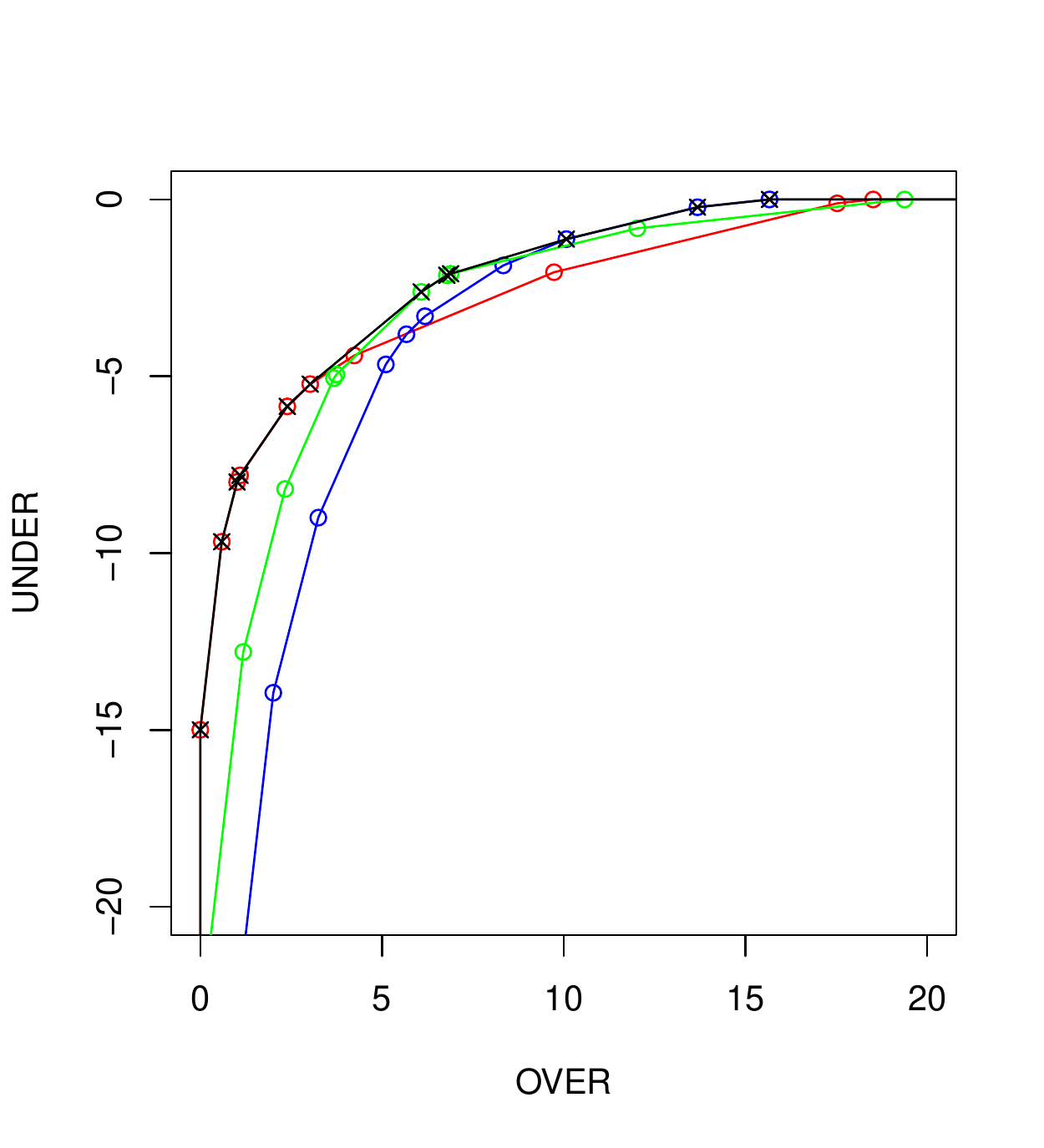} 
\vspace{-0.5cm}
\caption{Convex hull of Figure \ref{fig:RROCshiftall}, shown in black. There are 12 visible points (represented as black crosses) on the convex hull: 6 from $m_1$ (red), 3 from $m_3$ (green), and 3 from $m_2$ (blue).}
\label{fig:RROC3CH}
\end{figure}

From these three curves, we can calculate their convex hull, as shown in Figure \ref{fig:RROC3CH}.

\subsection{Areas and metrics}

\RROC analysis, as ROC analysis, can be especially useful for analysing models under different operating conditions and select the best one for a single operating condition or a region, or even better, to create hybrids through the notion of convex hull. Nonetheless, in ROC analysis we are also interested in evaluating models that can work well for a wide range of operating conditions. One measure that gives us a good indication of a classifier performing well in a wide range of operating conditions is the Area Under the ROC Curve (\auc). Can we develop a similar measure for \RROC curves?

The good mapping so far between ROC curves and \RROC curves in terms of what they represent suggests that this is possible. The following definition introduces such a measure:

\begin{definition}
The Area Over the \RROC Curve (\aoc) is defined as follows:
\[ \aoc \triangleq - \int_0^\infty \UNDER \:\:\: d\OVER = \int_{-\infty}^0 \OVER \:\:\: d\UNDER\]
\end{definition}

Lower values for \aoc are better.

The previous area can be calculated very easily using the sum of the $n+1$ upward trapeziums given between the elements 1 and n+2 from $RROCX$ and $RROCX$ in algorithm \ref{algo:RROCcurve}. Actually, for models always outputting finite values, this can be calculated from 2 to $n$, since the extreme trapezium 1 to 2 has area 0 and the trapezium $n+1$ to $n+2$ as well, so this only need to sum $n-1$ trapeziums. Consequently:

\[ \aoc = \sum_{i=2}^{n}-\frac{RROCY_{i+1}+RROCY_{i}}{2}(RROCX_{i+1}-RROCX_{i}) \]

The first question about this area is why we have defined the area over the curve and not under the curve. This has an easy answer: since the \RROC space is unbounded, the area under the curve is always infinite. But what about the \aoc? The following proposition gives an answer:

\begin{proposition}
For any regression model $m$ which always outputs finite values, the \aoc is finite.
\end{proposition}
\begin{proof}
Since the model $m$ always outputs finite values, there is a shift $s_o$, such that for any shift $s\leq s_o$ we have that $\OVER=0$ and there is also a shift $s_u$, such that for any shift $s\geq s_u$ we have that $\UNDER=0$. This means that the curve touches (and stays at) both the \xaxis and the \yaxis. Then the area is finite.
\end{proof}

For the three models in Figure \ref{fig:RROCshiftall}, the \aoc is 56.1387, 88.0933 and 63.9295 for models $m_1$, $m_2$ and $m_3$ respectively. Although a single number loses most of the information we can see on the curve, these numbers summarise their overall performance.

From the notion of \aoc, we can investigate what exactly means to have low \aoc and high \aoc.
The `best' model in terms of \aoc (a perfect square with top-left corner at the \RROC heaven (0,0)) means that there is a shift that achieves 0 error. This is rarely the case, except for datasets for one single example (where there is always a shift getting 0 loss). It is also very rare to have a dataset for which the error is always the same, another possible situation where we would have 0 \aoc.
 Note that a model with very high \mse or \mae could, in principle, have $\aoc=0$. This would suggest that the shift was very badly chosen. The parallel with classical ROC analysis here is clear, where we can have bad accuracy for a model with optimal \auc by choosing a bad threshold.

\begin{figure}
\centering
\includegraphics[width=0.5\textwidth]{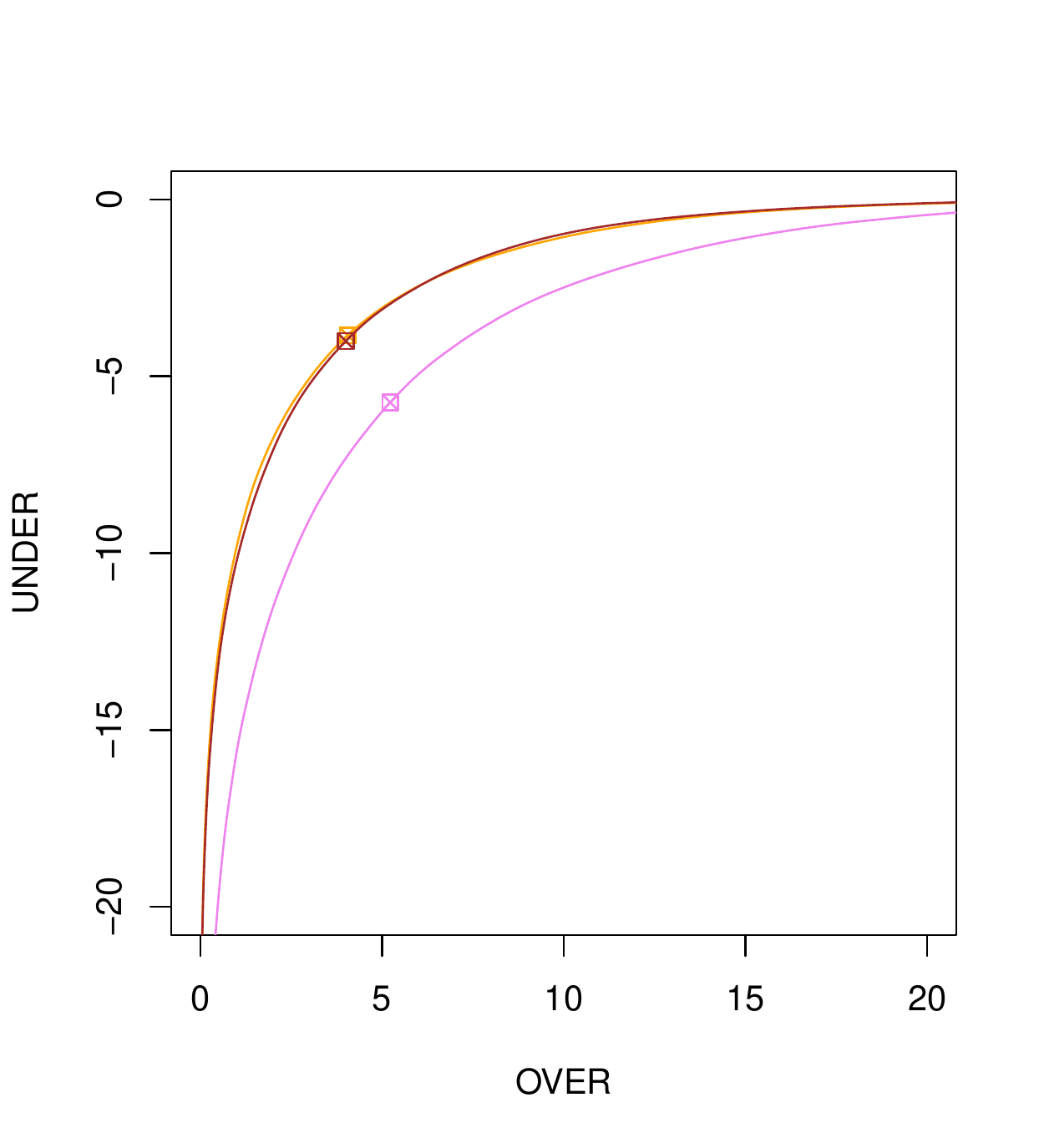} 
\vspace{-0.5cm}
\caption{We use a dataset with $1,000$ examples using a normal distribution with mean at 0 and standard deviation $0.01$. The first model (in violet) is just generated as a random model using the same distribution ($\mae= 0.011$, $\mse= 0.00019$, $\aoc=93.28$). As we see, its performance is low. The second model (in orange) is just built by using the actual values plus a random normal noise with the same distribution ($\mae= 0.0079$, $\mse= 0.00010$, $\aoc=49.83$). Finally, the third model (in brown) is just a model always outputting the same value, the mean of the actual data ($\mae= 0.0080$, $\mse=0.00010$, $\aoc= 50.31$). These two last models seem to have very similar metrics and curves (they almost overlap completely).}
\label{fig:RROCnormal}
\end{figure}

What about the `worst' model in terms of \aoc? Of course we can have a value of \aoc as high as we want. We can even get an infinite \aoc, if the model outputs $\infty$ or $-\infty$ for one single example. So, the question must be stated more precisely: given a model with a certain \mae, what is the worst value for \aoc?
This is difficult to answer. At first sight, it seems that the degree of dispersion of the error may affect, since it may make the shift more effective. Also, the degree of correlation between the actual and predicted values could be important.
Figure \ref{fig:RROCnormal} shows how a random model looks (in violet), which typically shows low performance. Also, it compares two models with similar performance, but one which is just generated adding random noise to the true values (in orange) and the other by calculating the mean of the true values (in brown). While these two last models have very different dispersion (the last model has null dispersion) and very different correlation (the last model has null correlation), their metrics and \RROC curves are very similar. This is explained because their error distributions are similar. Hence, one possible way of looking at \RROC curves is precisely this. They represent the distribution of errors.

A different question is to give a numerical interpretation of the \aoc. While its definition suggests that it may be the expected value of the total under-estimation given a uniform value for the total over-estimation, this is not well-defined because both \UNDER and \OVER are not bounded. A possible interpretation is that it aggregates the macro-average squared error ($n \cdot\mmse$) with a distribution which depends on the model\footnote{Note that the length of each segment may represent the frequency of each possible value of the asymmetry parameter $\alpha$.}, which is similar to one recent interpretation given to \auc \cite{hand2009measuring}.  
Other interpretations as an aggregation of expected loss may be possible\footnote{We suggest some possible pathways for exploration. Since \aoc is related to the magnitude of predictions (and errors), it cannot be directly related to rank-based correlation measures in regression. However, it could be related to this sum $\sum_{y_1,y_2}(\hat{y}_1 - \hat{y}_2  \:\:|\:\:   y_1 > y_2)$, which would work as a counterpart of the Wilcoxon-Mann-Whitney interpretation of the \auc.}, as it has happened to \auc recently, where new interpretations have been introduced 
\cite{ICML11CoherentAUC,JMLR12}.

Having said all this, our previous idea of the \aoc being related to the distribution of errors seems more appealing. If we have a compact error distribution, then \aoc will be low. If we have a sparse error distribution, then \aoc will be high. One classical measure of dispersion is precisely the variance, defined and decomposed as follows:

\begin{definition}\label{def:var}
The error variance $\sigma(\vect{e})^2$ is defined as:
\[ \sigma(\vect{e})^2 \triangleq \frac{\sum_{i=1}^n ({e_i-\mu(\vect{e})})^2} {n} = \frac{\sum_{i=1}^n e_i^2} {n} - \mu(\vect{e})^2 = \mu(\vect{e}^2) - \mu(\vect{e})^2  = \mse(\vect{e}) - \mu(\vect{e})^2 \]
where $\mu(\vect{e})$ represents the mean of the vector $\vect{e}$.
\end{definition}

Note that we define the population variance, by dividing by $n$ (instead of $n-1$). The reason is just to keep the expressions that will follow next as simple as possible.
We will use just $\sigma$ (instead of $\sigma(\vect{e})$) and $\mu$ (instead of $\mu(\vect{e})$) when clear from the context.
The last term in definition \ref{def:var} is just a different way of showing the classical \mse decomposition as the sum of the squared error bias ($\mu^2$) and the error variance ($\sigma^2$).

Quite surprisingly, the observation that the \aoc and the error variance are related can be made extremely precise, as the following theorem shows:

\begin{theorem}\label{theo:wow}
The area over the \RROC curve equals the population variance $\sigma^2$ of the errors multiplied by a factor $n^2/2$ which is independent of the model. Namely:
\[\aoc = \frac{\sigma^2 n^2}{2} \]
\end{theorem}
\begin{proof}

We start with an error vector $\vect{e}$ of length $n$, which we assume is sorted in decreasing order, as in algorithm \ref{algo:RROCcurve}.
We use a different notation for the points in the \RROC curve. Instead of using $n+2$ points, we will just ignore the two extremes (which do not contribute to the area for finite cases) and we will just work with $n$ points, denoted by $p_1, \dots, p_n$. The components of each point are $p_i = (o_i, u_i)$. Note that $o_i = RROCX_{i+1}$ using the notation in algorithm \ref{algo:RROCcurve} and $u_i = RROCY_{i+1}$.
We will also introduce the error differences $d_i = e_{i} - e_{i+1}$, which are defined from $i=1$ to $i=n-1$. Note that $d_i \geq 0$ since the error vector $\vect{e}$ is in decreasing order.
It is easy to see that $o_i = \sum_{j=1}^{i-1} j \cdot d_j$ and $u_i = -\sum_{j=i}^{n-1} (n-j)\cdot d_j$.
According to this notation:

\[ \aoc = - \sum_{i=1}^{n-1} \frac{u_i + u_{i+1}}{2} (o_{i+1} - o_{i}) \]

In order to prove this theorem, we will proceed by induction.

\noindent {\bf Base case}

The base case will consider any error vector of size $n=2$. In this case, we only have two points $p_1 = (0, -d_1)$ and $p_2 = (d_1, 0)$. From here,

\begin{eqnarray*}
\aoc & = & -\sum_{i=1}^{1} \frac{u_i + u_{i+1}}{2} (o_{i+1} - o_{i}) = - \frac{-d_1 + 0}{2} (d_1 - 0) = \frac{d_1^2}{2} = \frac{(e_2 - e_1)^2}{2}  \\
     & = &  =	\frac{(e_2 - \mu + \mu - e_1)^2}{2} = \frac{ (e_2 - \mu)^2 + (\mu - e_1)^2 + 2(e_2 - \mu)(\mu - e_1) }{2}\\
		 & = &  = \frac{2(e_2 - \mu)^2 + 2(\mu - e_1)^2}{2} = \frac{4\sigma^2}{2} = \frac{\sigma^2 n^2}{2} 
\end{eqnarray*}


\noindent {\bf Inductive step}

We assume that 
\begin{eqnarray}
\aoc = \frac{\sigma^2 n^2}{2} \label{eq:inductivestep}
\end{eqnarray}
holds for any dataset of size $n$.

Without loss of generality, we consider that the case for $n+1$ is constructed by adding example $e_{n+1}$, assumed to be lower than the other examples $e_1, e_2, \dots, e_n$ in the case for $n$. Consequently, the error vector for the case $n+1$ is $e_1, e_2, \dots, e_n, e_{n+1}$. The difference vector is also an extension for $n+1$, denoted by $d_1, d_2, \dots, d_n$.
Note that since we assume that eq. \ref{eq:inductivestep} holds for any dataset of $n$ examples, we can choose the order of examples that we prefer in order to build any case with $n+1$ examples.

The $\aoc$ for the $n$ case is given by:
\[ \aoc = - \sum_{i=1}^{n-1} \frac{u_i + u_{i+1}}{2} (o_{i+1} - o_{i}) \]

The $\aoc$ for the $n+1$ case is given by
\begin{equation}
\widetilde{\aoc} = - \sum_{i=1}^{n} \frac{\widetilde{u}_i + \widetilde{u}_{i+1}}{2} (\widetilde{o}_{i+1} - \widetilde{o}_{i}) \label{eq:newaoc}
\end{equation}

We will use a wide tilde to denote the $\widetilde{\aoc}$, $\widetilde{\sigma}$, $\widetilde{\mu}$, etc., for the $n+1$ case.
The first thing we can see is that $\widetilde{u}_1 = u_1 - d_1 - d_2 - \dots - d_n$, $\widetilde{u}_2 = u_2 - d_2 - \dots - d_n$, etc.
We use these latter expressions on (\ref{eq:newaoc}):

\[
\widetilde{\aoc} = - \sum_{i=1}^{n} \frac{u_i + \{ -\sum_{j=i}^n d_j \} + u_{i+1} + \{ -\sum_{j=i+1}^n d_j\}}{2} (\widetilde{o}_{i+1} - \widetilde{o}_{i}) 
\]

The second thing we realise is that $o_i$ and $\widetilde{o}_i$ are equal for $i=1 \dots n$. From here, we can calculate the delta between $n+1$ and $n$ as follows:

\[
\Delta\aoc \triangleq \widetilde{\aoc}- \aoc = - \sum_{i=1}^{n} \frac{\{ -\sum_{j=i}^n d_j \} + \{ -\sum_{j=i+1}^n d_j\}}{2} (\widetilde{o}_{i+1} - \widetilde{o}_{i}) 
\]

But we have that $\widetilde{o}_{i+1} - \widetilde{o}_{i} = i \cdot d_i$. So, we rewrite:

\begin{eqnarray*}
\Delta\aoc & = & - \sum_{i=1}^{n} \frac{\{ -\sum_{j=i}^n d_j \} + \{ -\sum_{j=i+1}^n d_j\}}{2} (i \cdot d_i) \\
           & = & \sum_{i=1}^{n} \frac{ (d_i + 2 \sum_{j={i+1}}^n d_j) (i \cdot d_i)  }{2}  \\
					 & = & \sum_{i=1}^{n} \frac{ i \cdot d_i^2 + 2 \sum_{j={i+1}}^n i \cdot d_i d_j  }{2} 
 \end{eqnarray*}

Using the expression of the square of a sum: $(\sum_i a_i)^2 = \sum_i a_i^2 + 2 \sum_{i < j} a_i a_j$, and joining/distributing terms, we see that the above expression can be rewritten as:

\begin{eqnarray*}
\Delta\aoc & = & \sum_{i=1}^{n} \frac{ \{\sum_{j=i}^n d_i \}^2  }{2} \\
           & = & \sum_{i=1}^{n} \frac{ (e_{n+1} - e_i)^2  }{2} \\
           & = & \sum_{i=1}^{n} \frac{ e_{n+1}^2 - 2e_{n+1} e_i + e_i^2  }{2} \\
           & = & \frac{1}{2} \left(n \cdot e_{n+1}^2 - 2e_{n+1} \sum_{i=1}^{n} e_i + \sum_{i=1}^{n} e_i^2 \right)  \\
           & = & \frac{1}{2} \left(n \cdot e_{n+1}^2 - 2e_{n+1} n \cdot \mu    + n(\sigma^2+\mu^2)  \right)  \\					
           & = & \frac{n}{2} \left(e_{n+1}^2 - 2e_{n+1} \mu    + \sigma^2+\mu^2  \right)  \\
           & = & \frac{n}{2} \left( (e_{n+1} - \mu)^2    + \sigma^2  \right)  
 \end{eqnarray*}

From here, we can now write:

\begin{eqnarray*}
\widetilde{\aoc} = \aoc + \Delta\aoc = \aoc + \frac{n}{2} \left( (e_{n+1} - \mu)^2    + \sigma^2  \right)  
\end{eqnarray*}

From the induction step (equation \ref{eq:inductivestep}), we have:

\begin{eqnarray*}
\widetilde{\aoc} & = & \frac{\sigma^2 n^2}{2} + \frac{n}{2} \left( (e_{n+1} - \mu)^2    + \sigma^2  \right)  \\ 
                 & = & \frac{n}{2} \left( \sigma^2 n + (e_{n+1} - \mu)^2    + \sigma^2  \right)  \\ 
                 & = & \frac{n}{2} \left( \sigma^2 (n+1) + (e_{n+1} - \mu)^2      \right)  \\ 
                 & = & \frac{n}{2} \left( (\frac{\sum_{i=1}^n e_i^2}{n}-\mu^2) (n+1) + (e_{n+1} - \mu)^2      \right)  \\ 
                 & = & \frac{1}{2} \left( ({\sum_{i=1}^n e_i^2}-n\mu^2) (n+1) + n\cdot (e_{n+1} - \mu)^2      \right)  \\ 
								& = & \frac{1}{2} \left( ({\sum_{i=1}^n e_i^2})(n+1)- (n\mu^2) (n+1) + n\cdot e_{n+1}^2 - 2n \cdot e_{n+1} \mu + n\mu^2      \right)  \\ 
								& = & \frac{1}{2} \left( ({\sum_{i=1}^{n+1} e_i^2})(n+1)- (n\mu^2) (n) - e_{n+1}^2 - 2n \cdot e_{n+1} \mu       \right)  \\ 
								& = & \frac{1}{2} \left( ({\sum_{i=1}^{n+1} e_i^2})(n+1) - (n\mu  + e_{n+1})^2 \right)  \\ 
								& = & \frac{1}{2} \left( ({\sum_{i=1}^{n+1} e_i^2})(n+1) - ((n+1) \widetilde{\mu})^2 \right)  \\
								& = & \frac{1}{2} (n+1)^2 \left( \frac{{\sum_{i=1}^{n+1} e_i^2}}{n+1} - \widetilde{\mu}^2 \right)  \\
								& = & \frac{1}{2} (n+1)^2 \left( \widetilde{\sigma}^2 \right)  \\
								& = & \frac{\widetilde{\sigma}^2 (n+1)^2}{2}
\end{eqnarray*}
 
This last expression completes the induction step and so does the proof.
\end{proof}

\begin{corollary}\label{cor:unbiased}
If the model is unbiased (i.e. $\mu(\vect{e})=0$) then:
\[\aoc = \frac{\mse \cdot n^2}{2} \]
\end{corollary}
\begin{proof}
The proof is direct from theorem \ref{theo:wow} and definition \ref{def:var} (\mse decomposition).
\end{proof}

For the models $m_1$, $m_2$ and $m_3$ in examples \ref{ex:example1}, \ref{ex:example2} and \ref{ex:example3} we have a variance of $1.1228$, $1.7619$ and $1.2786$ respectively. The \aoc is $56.1387$, $88.0933$ and $63.9295$ respectively. Since $m_2$ is unbiased, its $\mse$ is precisely $1.7619$, its error variance. The constant factor is $n^2/2=50$ in the three cases.
Similarly, for the third model (in brown) in Figure \ref{fig:RROCnormal}, as it always outputs the mean of a distribution with standard deviation $0.01$ and $1,000$ examples, we have that the \aoc was 50.31. The expected result is $0.001 \cdot 1000^2/2= 50$. The difference is not given because theorem \ref{theo:wow} is approximate, but just because the sample is generated with a distribution with $\sigma^2= 0.0001$, but the sample does not exactly have this variance (it is actually $0.00010062$).

\begin{figure}
\centering
\includegraphics[width=0.5\textwidth]{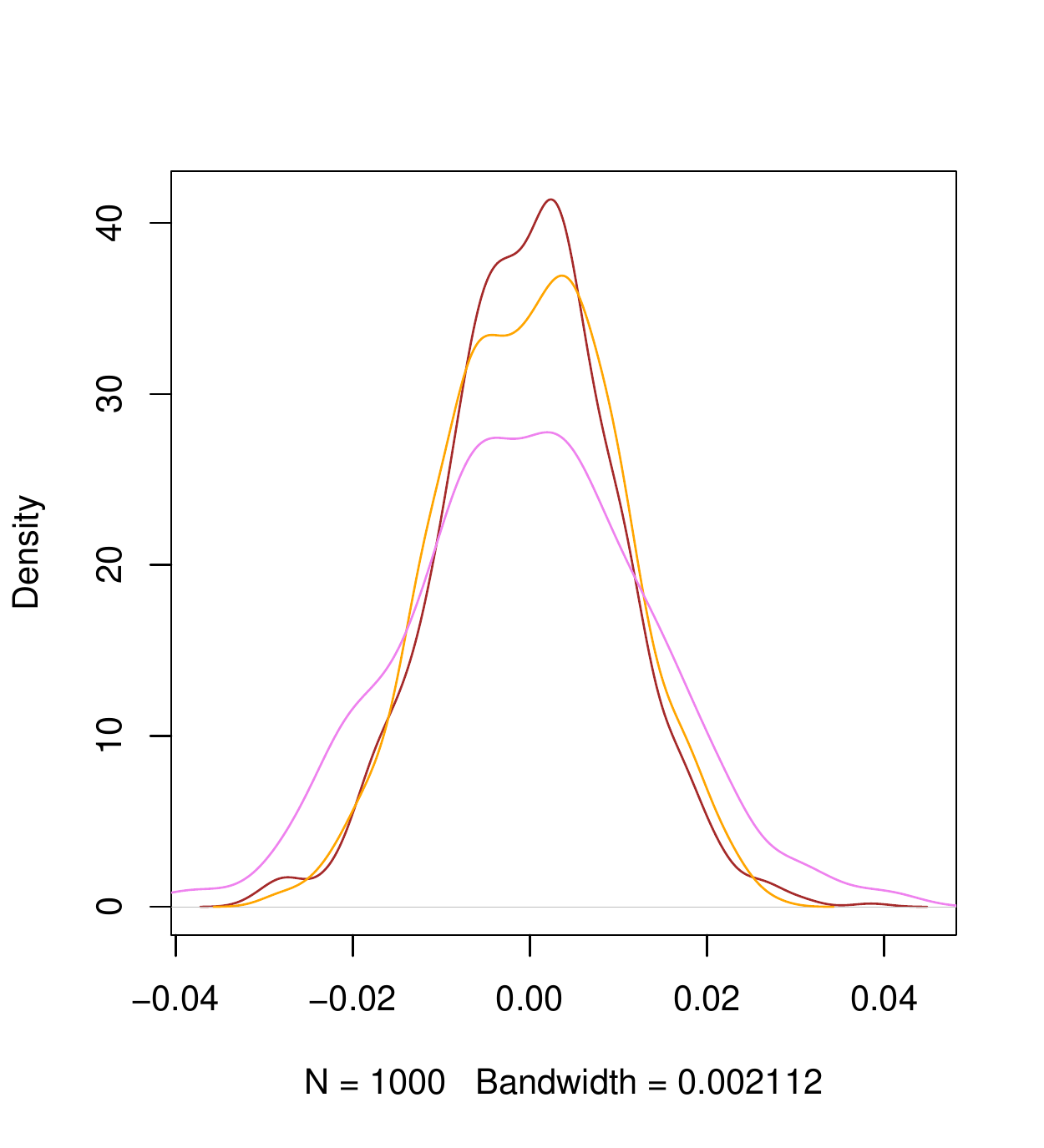} 
\vspace{-0.5cm}
\caption{The error density plots of the same three models as in Figure \ref{fig:RROCshiftall}. Note that we can see that similar error density functions will produce similar \RROC curves. Here, more peaked density functions mean better \RROC curves.}
\label{fig:density}
\end{figure}

Given the connection between the area over the \RROC curve and the population variance, we can explore the connection between the \RROC curve and an error density plot. As we can see in Figure \ref{fig:density}, there is a high correspondence between the density plots and the \RROC curve, but the cumulative character of the \RROC curve make the latter smoother.

Note that this connection between $\aoc$ and the error variance indicates that it is the dispersion that counts when trying to adapt our models to cost-sensitive situations with asymmetries, and not the position, which can be ignored by assuming that the optimal shift will be chosen for each particular operating condition. This again shows a parallel with ROC analysis. In ROC analysis, the absolute values of the scores do not affect the \auc. Only their order matters. Here, for \RROC curves, the position of the mean error (the error bias) does not affect the \aoc, only the dispersion of the error.

\begin{figure}
\centering
\includegraphics[width=0.5\textwidth]{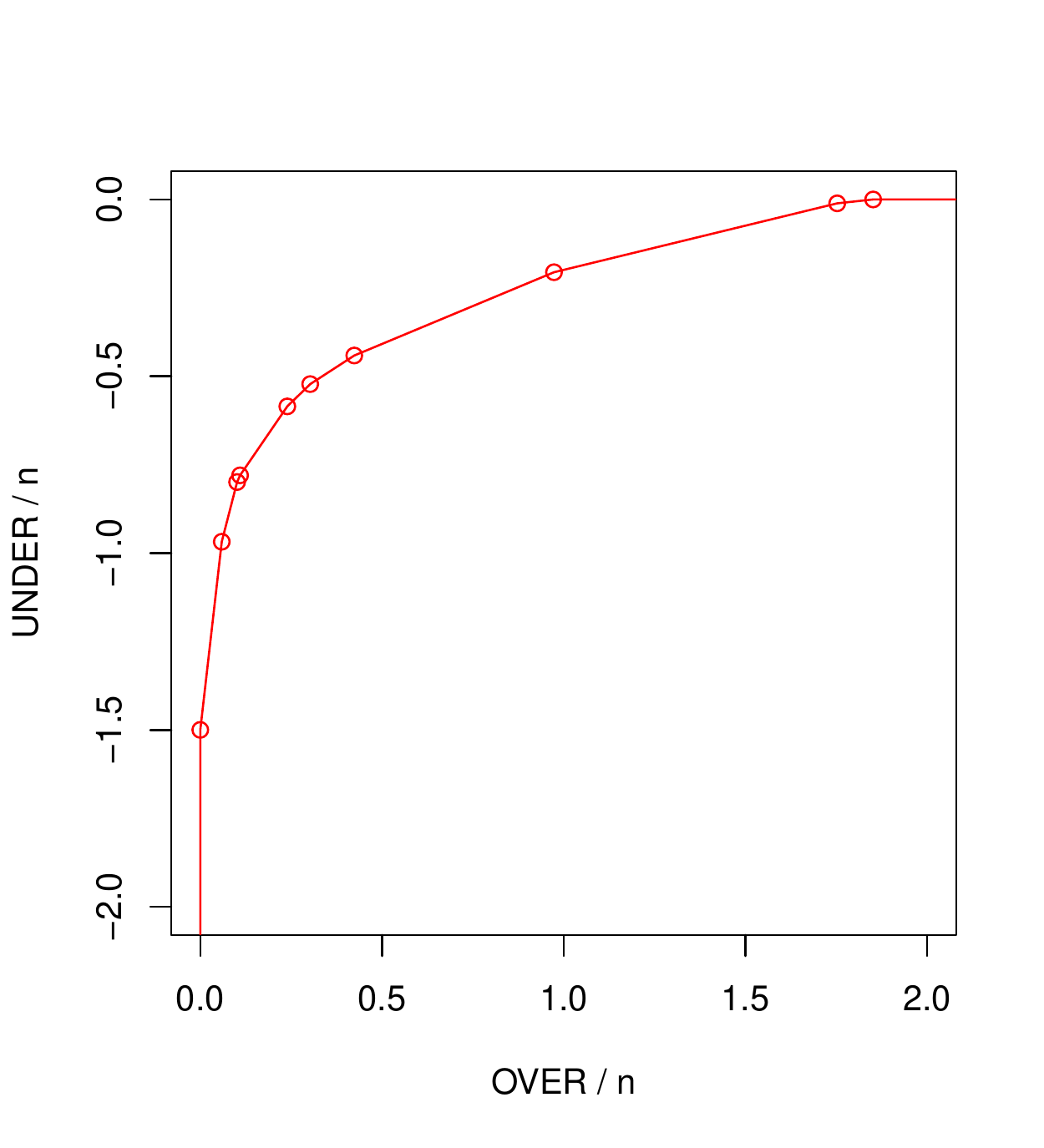} 
\vspace{-0.5cm}
\caption{Model $m_1$ as in example \ref{ex:example1} drawn as a normalised \RROC curve, using a normalised scale for the \xaxis and \yaxis (n=$10$ examples).}
\label{fig:RROCnormalised1}
\end{figure}

This is a fundamental result as well because it is a graphical representation of the error variance, which can sum up to the applicability of \RROC curves. The $n^2$ factor in theorem \ref{theo:wow} also suggests that a scaled representation of \RROC curves could be done by dividing both the \xaxis and \yaxis by $n$, i.e., plotting $\OVER/n$ against $\UNDER/n$. This would make the curves independent of the number of examples, but the meaning of each point would be somewhat blurred, as the `average over-estimation or (under-estimation) per example'. Nonetheless, this could be the standard representation in many application, especially when the number of examples in the datasets may vary or we may even compare several models (or the same model) against different datasets (with different sizes).
Figure \ref{fig:RROCnormalised1} shows the same plot as Figure \ref{fig:RROCshift1} but normalising by the number of examples (in this case $n=10$).

\section{ROC curves for non-constant shifts and soft regressors}\label{sec:nonconstant}

\begin{figure}
\centering
\includegraphics[width=0.5\textwidth]{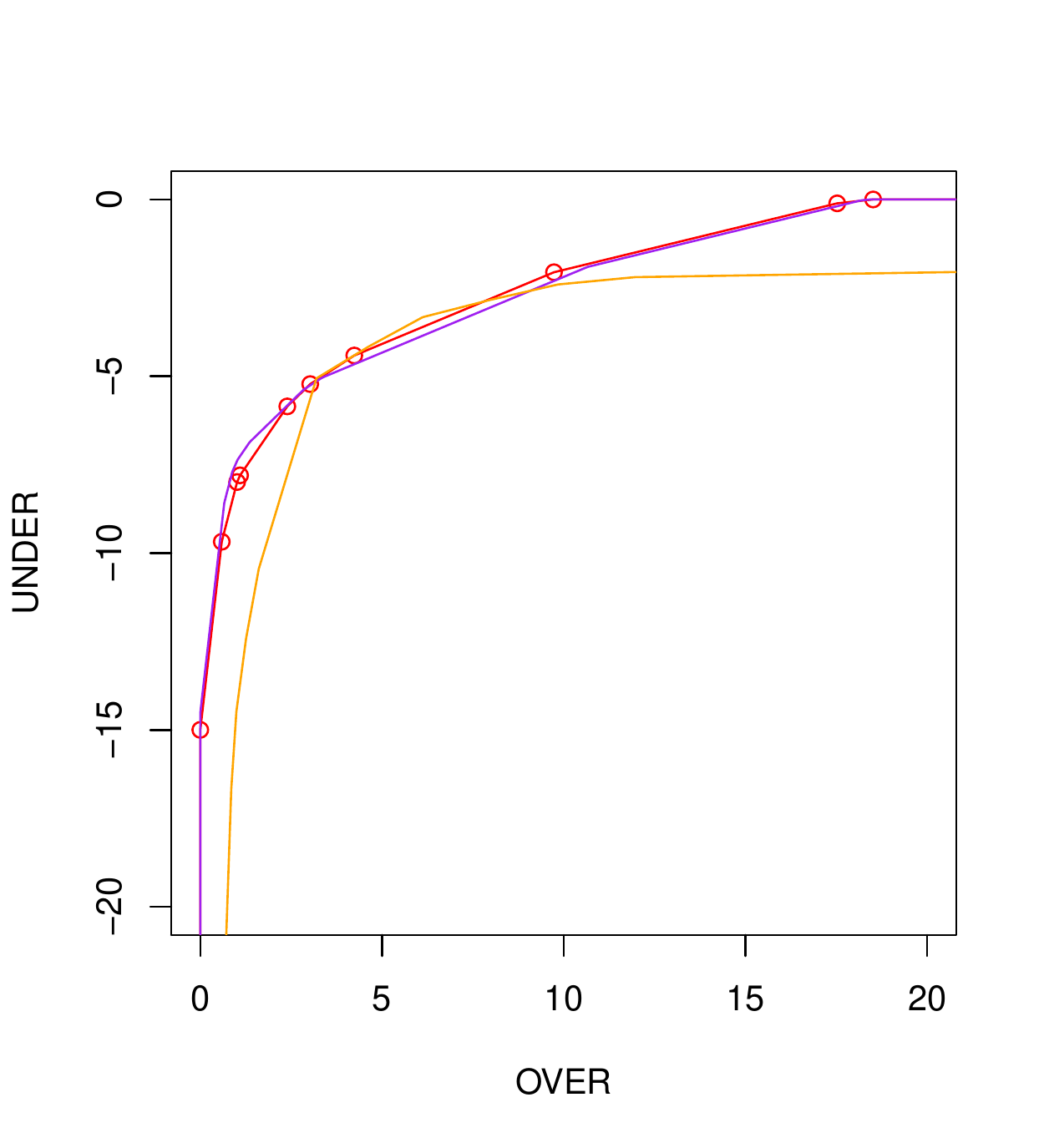} 
\vspace{-0.5cm}
\caption{\RROC Curve of model $m_1$ from example \ref{ex:example1} using a constant shift as usual: $\hat{y} + s$ (red), the same model using a third-degree polynomial, as $0.9\cdot\hat{y} + 0.002\cdot\hat{y}^3 + s$  (purple), and the same model using a third-degree polynomial combining $s$ and $\hat{y}$, as $\hat{y} - 0.004\cdot\hat{y}^3 + 1.5\cdot s \cdot\hat{y}^2 + s$  (orange).}
\label{fig:RROCnonconstant}
\end{figure}

In classification, there are many possibilities for choosing the threshold \cite{JMLR12}.
In regression, there are many possibilities as well for the shift. Until now, we have considered that the shift is chosen as a constant. Other possibilities rely on the use of any function of the prediction and the operating condition. Figure \ref{fig:RROCnonconstant} shows the model $m_1$ from example \ref{ex:example1} using a constant shift, the same model using a third-degree polynomial, and the same model using a third-degree polynomial combining $s$ and $\hat{y}$. 
As we can see on the figure, there are places where the use of a different shift formula can reach places where the constant shift cannot.
Actually, we can find functions such that the predictions are modified in such a way that they can attain any point on the \RROC space. However, in order to get close to the \RROC heaven, we would need very ad-hoc functions, basically embedding an error correction inside.

In general, we are interested in shift functions and methods that are systematic (i.e., a procedure which is the same for all models). Clearly, a constant shift is a systematic method, provided we find a way to find the appropriate constant for each operating condition. Recently, a method to find the appropriate shift for each operating condition (asymmetry) has been introduced \cite{Bansal2008}. Simply, given a value of $\alpha$, the method calculates the best shift for the training set. Then, this shift can be applied to the test set. This method does not obtain the optimal shift for every $\alpha$, but if the training set and test set are similar, the approximation can be good. By ranging over operating conditions $\alpha$ (instead of shifts), and using this method, we can construct a \RROC curve which does not show the evolution of \OVER and \UNDER for an optimal (or ideal) shift choice method, but an actual, feasible one. This reinforces the view expressed by \cite{JMLR12} for classification: we evaluate pairs of models and threshold choice methods. The translation to regression and \RROC curves is that we plot models assuming a shift choice method (both threshold choice methods and shift choice methods are types of reframing methods). In the previous sections, we have assumed an optimal constant-shift choice method, but many other options exist and may lead to other curves for the same model (which are not necessarily convex). The good thing about the \RROC space is that we can visualise several options in the same plot, as done with Figure \ref{fig:RROCnonconstant}, and evaluate both models and shift choice methods at the same time.

Overall, there are many shift choice methods to be explored. For instance, \cite{zhao2011extended} generalise the constant-shift choice method from  \cite{Bansal2008} by using any polynomial function. A different, and more powerful, perspective is introduced by \cite{probreg2012}, where instead of using crisp models, the regression model accompanies a standard deviation to each prediction. This standard deviation is used to better adjust the shift according to each example, which is now a function of two variables instead of one. The adjustment is found by  risk minimisation.

This also suggests the exploration of the connection between \RROC curves and its corresponding cost curves.

\begin{definition}
The cost space for regression is defined as a plot where the expected loss (e.g., the asymmetric absolute loss) is shown on the \yaxis for a range of operating conditions (e.g., the asymmetry $\alpha$).
\end{definition}

\begin{figure}
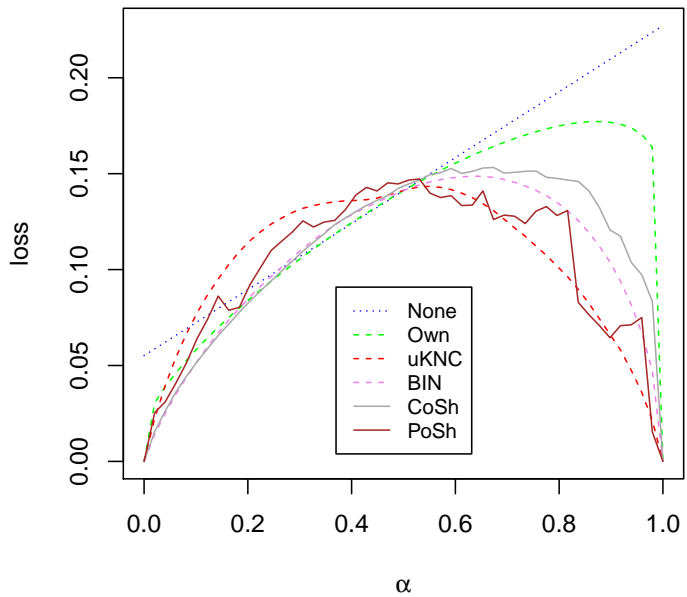

\centering
\includegraphics[width=0.6\textwidth]{{{D01-iris3-LR-SDMET3.1WEXP0-CPAg0.5r0.5-ABSLOSS-alphas50rhos1-normal-noshuffleCV-RO2}}} 
\vspace{-0.5cm}
\caption{A cost curve (\cite{probreg2012}) showing the absolute loss $\aabsloss$ (definition \ref{def:absloss}) against the operating condition $\alpha$ (asymmetry) of a model using several shift choice methods (reframing methods). The `None' method represents no adjustment and corresponds to a single point in \RROC space. The `CoSh' and `PoSh' are the constant-shift choice method from \cite{Bansal2008} and the polynomial-shift choice method from \cite{zhao2011extended}. Finally, the `Own', `uKNC' and `BIN' methods are probabilistic reframing methods based on soft regressors using a two-parameter output for each instance and a risk-minimisation solution. Note that the optimal constant-shift choice method is not shown here, but it should fall below the `CoSh' and `PoSh' methods. }
\label{fig:costcurve}
\end{figure}

Figure \ref{fig:costcurve} shows this cost space, which is similar to the cost space of Drummond and Holte's cost curves for classification\cite{drummond-and-Holte2006}. The investigation of the mapping between the regression cost space and the \RROC space can lead to new important findings as has been recently done for classification \cite{ROCandCost}.

\section{Concluding remarks}\label{sec:conclusions}

We said in the introduction that there is no such a thing as the `canonical' ROC space for regression, corresponding exactly to the ROC space for classification, since regression and classification are very different tasks. Having said this, we think that the \RROC space, curves and analysis that we have introduced in this paper present so many parallelisms and share so many notions and procedures, that their curves could reasonably called {\em the} ROC curves for regression, with arguable more support than other previous attempts. We have seen that the notions of operating condition, cost asymmetry, \RROC space, points, segments, \RROC heaven, \RROC isometrics, hybrid models, convexity, dominance, convex hull, curves, shift choice methods, etc., derive smoothly and work almost the same as in the classification case, so the practitioners which are used to ROC curves can directly apply their expertise on ROC analyse to regression quite easily.

There are naturally several issues which could lead to more general (or slightly different) notions of \RROC curve for regression, keeping the same basic structure. The first issue that could be explored and generalised is the very notion of operating condition. We have only considered the asymmetry while, in classification, the class distribution can also be integrated (along with the cost proportion) in what is usually referred to as {\em skew}.
In regression, the distribution of the output value (and not only the loss asymmetry) may also be considered part of the operating condition as well. This integration does not seem to be direct, but it is worth being investigated.

A second issue is the use of other loss functions. For instance, instead of an asymmetric absolute error, we could use an asymmetric squared error {\em Quad-Quad}. We guess that this would lead to non-straight isometrics and non-straight segments in the \RROC curve, but the basic ideas would remain. Again, plotting different isometrics in \RROC space for many different loss functions ({\em Lin-Lin, Quad-Quad, Lin-Exp, Quad-Exp, etc.}) would be a work on its own, very much resembling the 
celebrated paper \cite{Fla03} on isometrics for ROC curves in classification.



A third important avenue of future work would be to further investigate the connection with the error variance we have unveiled here and to analyse the relation of \aoc to other metrics, as well as the relation of \RROC space with other plots to analyse the performance of regression. We think that \RROC curves represent the expected loss for a range of operating conditions on one side, and the distribution of the error on the other side. There may be important connections to be unveiled between regression techniques trying to minimise the error variance (which we have shown here to be equal to the AOC) instead of squared error and those classification techniques trying to maximise the AUC (which has recently been shown to be equivalent to the refinement loss term of the MSE decomposition using the ROC curve \cite{JMLR12}) instead of accuracy \cite{ferri2002learning}\arxiv{\cite{FFH03}}.
So we 
anticipate 
a plethora of connections between \RROC curves and many other performance metrics in regression, as has been done for classification in the past years \cite{PRL09,hand2009measuring,ICML11CoherentAUC,ICML11Brier,JMLR12}.

Overall, we think that \RROC curves could become a fundamental tool in the assessment, improvement and deployment of regression models. In order to facilitate their use in real applications, we have developed a library for plotting \RROC curves, calculating their areas and deriving their convex hulls. The software, in R \cite{Rproject}, is available at \url{http://users.dsic.upv.es/~jorallo/RROC/}. The availability of software, the ubiquitous appearance of asymmetric losses in regression applications, and the success of ROC analysis for classification in the past decades suggests that \RROC curves may soon become mainstream in all the areas where ROC analysis has shown to be useful: medicine, bioinformatics, decision making, statistics, machine learning and pattern recognition.

\section*{Acknowledgments} 
I would like to thank Peter Flach and Nicolas Lachiche for some very useful comments and corrections on earlier versions of this paper, especially the suggestion of drawing normalised curves (dividing \xaxis and \yaxis by $n$).
This work was supported by the MEC/MINECO projects CONSOLIDER-INGENIO CSD2007-00022 and TIN 2010-21062-C02-02, GVA project PROMETEO/2008/051, the COST - European Cooperation in the field of Scientific and Technical Research IC0801 AT, and the REFRAME project granted by the European Coordinated Research on Long-term Challenges in Information and Communication Sciences \& Technologies ERA-Net (CHIST-ERA), and funded by the respective national research councils and ministries.

\bibliographystyle{plain}

\bibliography{biblio}


 \end{document}